\documentclass{amsart}          

\usepackage{amsmath,amsthm}     
\usepackage{amssymb}            
\usepackage{euscript}           
\usepackage{enumerate,calc}
\usepackage[matrix,arrow,curve,frame]{xy}    

\xymatrixcolsep{1.9pc}                          
\xymatrixrowsep{1.9pc}
\newdir{ >}{{}*!/-5pt/\dir{>}}                  

\newtheorem{thm}[subsection]{Theorem}
\newtheorem{defn}[subsection]{Definition}
\newtheorem{prop}[subsection]{Proposition}

\newtheorem{cor}[subsection]{Corollary}
\newtheorem{lemma}[subsection]{Lemma}

\theoremstyle{definition}  
\newtheorem{example}[subsection]{Example}

\newtheorem{remark}[subsection]{Remark}

  {\end{list}}

\newcommand{\dfn}{\textbf} 

\newcommand{\mdfn}[1]{\dfn{\mathversion{bold}#1}} 

\newcommand{\Smash}             {\wedge}

\newcommand{\tens}              {\otimes}               
\newcommand{\iso}               {\cong}  

\newcommand{\cat}{\EuScript}    
\newcommand{\cA}{{\cat A}}      
\newcommand{\cC}{{\cat C}}

\newcommand{\cE}{{\cat E}}

\newcommand{\Top}{{\cat Top}}
\newcommand{\Gtop}{\cat Top_G}
\newcommand{\Set}{{\cat Set}}

\newcommand{\sSet}{s{\cat Set}}

\newcommand{\field}[1]  {\mathbb #1} 
\newcommand{\A}         {\field A}
\newcommand{\F}         {\field F}
\newcommand{\R}         {\field R}

\newcommand{\Z}         {\field Z}
\newcommand{\C}         {\field C}
\newcommand{\Q}         {\field Q}

\DeclareMathOperator*{\colim}{colim}
\DeclareMathOperator*{\hocolim}{hocolim}
\DeclareMathOperator*{\holim}{holim}

\DeclareMathOperator{\coker}{coker}

\DeclareMathOperator{\uSigma}{\Sigma}

\newcommand{\ra}{\rightarrow}                   
\newcommand{\lra}{\longrightarrow}              
\newcommand{\la}{\leftarrow}                    
\newcommand{\lla}{\longleftarrow}               
\newcommand{\llra}[1]{\stackrel{#1}{\lra}}      
\newcommand{\llla}[1]{\stackrel{#1}{\lla}}      

\newcommand{\mapiso}{\llra{\cong}}
\newcommand{\mapisot}{\stackrel{\cong}{\rightarrow}}
\newcommand{\we}{\llra{\sim}}                   

\newcommand{\inc}{\hookrightarrow}              

\newcommand{\blank}{-}                          

\newcommand{\Cech}{\check{C}}
\newcommand{\CCech}{\v{C}ech\ }

\newcommand{\re}{Re}

\newcommand{\assign}{\mapsto}
\newcommand{\copr}{\,\amalg\,}

\newcommand{\he}{\simeq}

\newcommand{\pt}{pt}

\newcommand{\Zt}{{\Z/2}}
\newcommand{\Am}{\underline{A}}
\newcommand{\Zm}{\underline{\Z}}
\newcommand{\pim}{\underline{\pi}}
\newcommand{\Zc}{\underline{\Z}}
\newcommand{\Loop}{\Omega}
\newcommand{\SPi}{Sp^\infty}
\newcommand{\Spi}{Sp^\infty}
\newcommand{\SP}{SP}
\newcommand{\RP}{\R P}
\newcommand{\KK}{\Z\times BU}
\newcommand{\KKr}{\Z\times BU}

\newcommand{\upi}{\underline{\pi}}

\newcommand{\PK}{{\field P}}   
\newcommand{\PP}{P}             
\newcommand{\CP}{\C{\text P}}
\newcommand{\AG}{\text{AG}}
\newcommand{\KR}{KR}
\newcommand{\orb}{Or}

\newcommand{\tre}[1]{|#1|}

\newcommand{\Het}{H\Zc_{et}}
\newcommand{\HZ}{H\Zc}
\newcommand{\tH}{\tilde{H}}
\newcommand{\vf}{}
\numberwithin{equation}{section}

\begin{document}

\title{An Atiyah-Hirzebruch Spectral Sequence for $KR$-theory}

\date{April 6, 2003}
\author{Daniel Dugger}
\address{Department of Mathematics\\ University of Oregon\\ Eugene, OR
97403} 

\email{ddugger@math.uoregon.edu}


\maketitle

\tableofcontents

\section{Introduction}

In recent years much attention has been given to a certain spectral
sequence relating motivic cohomology to algebraic $K$-theory
\cite{Be,BL,FS,V3}.  This spectral sequence takes on the form
\[ H^p(X,\Z(-\tfrac{q}{2})) \Rightarrow K^{p+q}(X),\]
where the $H^s(X;\Z(t))$ are the bi-graded motivic cohomology groups, and
$K^n(X)$ denotes the algebraic $K$-theory of $X$.  It is useful in our
context to use topologists' notation and write $K^n(X)$ for what
$K$-theorists call $K_{-n}(X)$.  The above spectral sequence is the
analog of the classical Atiyah-Hirzebruch spectral sequence relating
ordinary singular cohomology to complex $K$-theory, in a way that is
explained further below.

It is well known that there are close similarities between motivic
homotopy theory and the equivariant homotopy theory of $\Zt$-spaces
(cf. \cite{HK1,HK2}, for example).
In fact there is even a forgetful map 
of the form
\[ (\text{motivic homotopy theory over $\R$}) \ra
(\text{$\Zt$-equivariant homotopy theory}), \]
discussed in \cite[Section 3.3]{MV} and \cite[Section 5]{DI}.  Our aim
in this paper is to construct the analog of the above motivic spectral
sequence in the $\Zt$-equivariant context. The spectral sequence takes
on the form
\[ H^{p,-\frac{q}{2}}(X,\Zm) \Rightarrow KR^{p+q}(X), \]
where the analog of algebraic $K$-theory is Atiyah's $KR$-theory
\cite{At}.  The analog of motivic cohomology is $RO(G)$-graded
Eilenberg-MacLane cohomology, with coefficients being the constant
Mackey functor $\Zm$.  The indexing conventions have been chosen for
their analogy with the motivic situation, and will be elucidated
further in just a moment.

The fact is that constructing the above spectral sequence is not at
all difficult, and there are many ways it could be done.  In
equivariant topology one has so many tools to work with that the
arguments end up being very simple.  Unfortunately most of these tools
are not yet available in the motivic context.  This paper tries to
develop the spectral sequence in a way that might eventually work
motivically, and which accentuates the basic properties of the
spectral sequence.  We introduce certain `twisted' Postnikov section
functors, and use these to construct a tower for the equivariant space
$\Z\times BU$ in which the layers are equivariant Eilenberg-MacLane
spaces.  The homotopy spectral sequence for the tower is essentially
what we're looking for---although technically speaking this only
produces half the spectral sequence, and to get the other half we must
stabilize.  The approach here is similar to the one advocated in
\cite{V2}, but was worked out independently (in fact there are several
differences, one being that \cite{V2} takes place in the stable
category).

\medskip

We'll now explain the methods of the paper in more detail, starting
with our basic notation.  Recall that every real vector space $V$ with
an involution gives rise to a $\Zt$-sphere $S^V$ by taking its
one-point compactification. If $\R$ and $\R_-$ denote the
one-dimensional vector spaces with trivial and sign involutions,
respectively, then any $V$ will decompose as $\R^p\oplus (\R_-)^q$ for
some $p$ and $q$.  So the spheres $S^V$ form a bi-graded family, and
when $V$ is as above we'll use the notation
\[ S^V = S^{p+q,q}. \]
Here the first index is the topological {\it dimension\/} of the
sphere, and the second index is called the {\it weight\/}.  Note that
when $V=\C^n$ (regarded as a real vector space with the conjugation
action) then $S^V = S^{2n,n}$; in particular, $\CP^1\cong S^\C =
S^{2,1}$.  The reader should be warned that this differs from the
bi-graded indexing introduced in \cite{At} and later used in
\cite{AM,Ar}.

Recall from \cite{LMS} that to give an equivariant spectrum $E$ is to
give an assignment $V\mapsto E_V$ together with suspension maps
$\Sigma^W E_V \ra E_{V\oplus W}$ which are compatible in a certain
sense.  One then has cohomology groups $E^V(X)$ for any representation
$V$, and when $V$ is as above we will likewise write
$E^V(X)=E^{p+q,q}(X)$ to correspond with our bi-graded indexing of the
spheres.  This `motivic indexing' is quite suggestive, and ends up
being a useful convention.

For the group $\Zt$, every real representation is contained in a
$\C^n$ for a large enough value of $n$.  The definition of equivariant
spectra can then be streamlined a bit by only giving the assignment
$\C^n \mapsto E_{(2n,n)}=E_n$ together with structure maps
$S^{2,1}\Smash E_n \ra E_{n+1}$.  This is the approach first taken in
\cite{AM}---albeit with different indexing conventions, as mentioned
above---and was later used in \cite{V1,J}.  We will treat spectra this
way throughout the paper.

\medskip

Our first example of such an object is the $KR$ spectrum.  The space
$\Z\times BU$ has an obvious $\Zt$-action coming from complex
conjugation on the unitary group $U$.  From another perspective, one
could model $\Z\times BU$ by the infinite complex Grassmannian, again
with the action of complex conjugation.  The reduced canonical line
bundle over $\CP^1$ is classified by an equivariant map $S^{2,1}=\CP^1
\ra \Z\times BU$, and so one gets $S^{2,1} \Smash (\Z\times BU) \ra
\Z \times BU$ by using the multiplication in $\Z\times BU$.  So we
have a $\Zt$-spectrum in which every term is $\Z\times BU$, and this
is called the $KR$ spectrum.  In fact, it is an Omega-spectrum:
equivariant Bott periodicity shows that the maps $\Z\times BU \ra
\Omega^{2,1}(\Z\times BU)$ are equivariant weak equivalences, or that
$KR^{s,t}(X) \cong KR^{s-2,t-1}(X)$.  The reference for this fact is
\cite{At}.

The second spectrum we will need is the equivariant Eilenberg-MacLane
spectrum $H\Zm$.  The easiest way to construct this, by analogy with
the non-equivariant case, is to consider the spectrum $\C^n \mapsto
\AG(S^{2n,n})$.  Here $\AG(X)$ denotes the free abelian group on the
space $X$, given a suitable topology.  The structure maps are the
obvious ones, induced in the end by the isomorphisms $S^{2,1}\Smash
S^{2n,n} \cong S^{2n+2,n+1}$.  It is proven in \cite{dS} that this
spectrum represents Eilenberg-MacLane cohomology with coefficients in
the constant Mackey functor $\Zm$, and that it is an Omega-spectrum.
The $n$th space $\AG(S^{2n,n})$ is therefore an equivariant
Eilenberg-MacLane space, and will be denoted $K(\Z(n),2n)$.  This is
the last of the basic notation needed to describe our results.

\medskip

Our goal in this paper will be to construct certain functors $P_{2n}$
on the category of $\Zt$-spaces, which are analogs of the classical
Postnikov section functors.  Roughly speaking, $P_{2n}X$ will be built
from $X$ by attaching cones on all maps from spheres `bigger than'
$S^{2n,n}$.  There are different possible choices for what is meant by
this, for which we refer the reader to Section 3.

As was pointed out above there is a Bott map $\beta\colon S^{2,1} \ra
\Z\times BU$ which classifies the reduced canonical line bundle over
$\CP^1$; let $\beta^n$ denote its $n$th power $S^{2n,n} \ra
\Z\times BU$.  Applying Postnikov section functors gives the induced
map $P_{2n}(S^{2n,n}) \ra P_{2n}(\Z\times BU)$.  The main goal of this
paper is the following:

\begin{thm}
\label{th:main}
There are Postnikov functors $P_{2n}$ on the category of $\Zt$-spaces
with the properties that
\begin{enumerate}[(a)]
\item $P_{2n}(S^{2n,n})$ is weakly equivalent to $K(\Z(n),2n)$, and
\item $P_{2n} (S^{2n,n}) \llra{\beta^n} P_{2n} (\Z\times BU) 
\ra P_{2n-2}(\Z\times BU)$
is a homotopy fiber sequence.
\end{enumerate}
\end{thm}

\begin{cor}
The tower
\[ \xymatrix{
 \cdots \ar[r] & 
P_4(\Z\times BU) \ar[r] & P_2(\Z\times BU) \ar[r] & P_0 (\Z\times BU)
}
\]
has the following properties:
\begin{enumerate}[(i)]
\item The homotopy fiber $F_n$ of the map $P_{2n}(\Z\times BU) \ra
P_{2n-2}(\Z\times BU)$ is an equivariant Eilenberg-MacLane space
$K(\Z(n),2n)$.  
\item The Adams operation $\psi^k\colon \Z\times BU \ra \Z\times BU$
induces a self-map of the tower, whose action on $F_n$ coincides with
the multiplication-by-$k^n$ map on the Eilenberg-MacLane space
$K(\Z(n),2n)$.
\end{enumerate}
\end{cor}

Looking at the homotopy spectral sequence of the above tower 
then gives the following:

\begin{cor}
There is a fringed spectral sequence $E_2^{p,q}\Rightarrow
KR^{p+q,0}(X)$ where $E_2^{p,q}= H^{p,-\frac{q}{2}}(X;\Zm)$ when
$p+q\leq 0$ and $q$ is even, and $E_2^{p,q}=0$ otherwise.  The
spectral sequence converges conditionally for $p+q<0$, is
multiplicative, and has an action of the Adams operations $\psi^k$ in
which $\psi^k$ acts on $E_2^{p,q}$ as multiplication by
$k^{-\frac{q}{2}}$.
\end{cor}

One can also stabilize the spectral sequence to avoid the awkward
truncation, but then one loses the action of the Adams operations.  To
this end, we let $W_n$ denote the homotopy fiber of $\Z\times BU \ra
P_{2n-2}(\Z\times BU)$.  In non-equivariant topology the $W_n$'s are
the connective covers of $\Z\times BU$, and are also the
spaces in the $\Omega$-spectrum for connective $K$-theory $bu$.  The
following result shows the same for the $\Zt$-case:

\begin{prop}
There are weak equivalences $W_n \ra \Loop^{2,1} W_{n+1}$, unique up
to homotopy, making the
diagrams
\[ \xymatrix{W_n \ar[r] \ar[d] & \Loop^{2,1}W_{n+1} \ar[d] \\
 \Z\times BU \ar[r] & \Loop^{2,1}(\Z\times BU) 
}
\]
commute (where the bottom map is the Bott periodicity map).
\end{prop}

The corresponding $\Zt$-spectrum whose $n$th object is $W_n$ will be
denoted \mdfn{$kr$} and called the \mdfn{connective $KR$-spectrum}.

\begin{thm}
\label{th:stabmain}
There is a `Bott map' $\beta\colon\Sigma^{2,1} kr \ra kr$ with the
following properties:
\begin{enumerate}[(a)]
\item The cofiber of $\beta$ is
$H\Zm$;
\item The telescope of the tower
\[ { \cdots \ra \Sigma^{2,1} kr \ra kr \ra \Sigma^{-2,-1} kr \ra
\Sigma^{-4,-2} kr \ra \cdots }
\]
is weakly equivalent to the spectrum $KR$ (where each map in the tower
is the obvious suspension or desuspension of $\beta$);
\item The homotopy inverse limit of the above tower is contractible.
\end{enumerate}
\end{thm}

The above tower of course yields a spectral sequence for computing
$KR^*(X)$ for any $\Zt$-space $X$, which could be considered the
Bockstein spectral sequence for the map $\beta$:

\begin{thm}
\label{th:stabspseq}
For any $\Zt$-space $X$, there is a conditionally convergent,
multiplicative spectral sequence
of the form $H^{p,-\frac{q}{2}}(X,\Zm) \Rightarrow KR^{p+q,0}(X)$.
\end{thm}

This spectral sequence is interesting even when $X$ is a point, in
which case it converges to the groups $KO^*$; it is drawn in detail in
section~\ref{se:spseq4pt}.  Also, note that there is really a whole
family of spectral sequences of the form
\[
H^{p,r-\frac{q}{2}}(X,\Zm) \Rightarrow KR^{p+q,r}(X), 
\]
but these can all be shifted back to the case $r=0$ by using Bott
periodicity $KR^{s,t}(X)=KR^{s+2,t+1}(X)$.

\subsection{Acknowledgments}
Most of the results in this paper were taken from the author's MIT
doctoral dissertation \cite{D1}.  The author would like to thank his
thesis advisor Mike Hopkins, and would also like to acknowledge very
helpful conversations with Gustavo Granja.  The final year of this
research was generously supported by a Sloan Dissertation Fellowship.

Since there has been a long delay between \cite{D1} and the appearance
of this paper, a brief history of related work might be in order.
Very shortly after \cite{D1} was written, Friedlander and Suslin
released \cite{FS} which constructed the more interesting motivic
spectral sequence, using very different methods.  In early 2000 the
paper \cite{V2} was released, outlining via conjectures a
homotopy-theoretic approach similar to the one given here (but working
in the stable category, and using a different definition of the
Postnikov sections).  These ideas were developed a little further in
\cite{V3}.  Sometime in 2000-2001 Hopkins and Morel also announced
proofs of results along these lines, although the details have yet to
appear.  At the end of 2002, the paper \cite{V4} proved a stable result
similar to Theorem~\ref{th:main}(a) in the motivic context, over
fields of characteristic zero.  An analog of Theorem~\ref{th:main} for
the unstable motivic category has never been claimed or proven, as far
as I know.

\subsection{Organization of the paper}
The paper has been written with a good deal of exposition, partly because
the literature on these subjects is not always so clear.  Sections 2 and 3
set down the necessary background, in particular giving the
constructions of equivariant Postnikov functors.  In these sections we
often work over an arbitrary finite group, because it is easier to
understand the ideas in this generality.  In this context everything
is graded by orthogonal $G$-representations, as is standard from
\cite{LMS}.  When specializing to the $\Z/2$ case we always translate
into the motivic $(p,q)$-indexing.  Section 2 also recalls the basic
facts we will need about the theory $H^{*,*}(\blank;\Zc)$.
 
The real work takes place in section 4, where we analyze the Postnikov
tower for $\Z\times BU$.  Section 5 discusses the basic properties of
the associated spectral sequence, most of which follow immediately
from the way the tower was constructed.  Section 6 is concerned with
passing to the stable case.  Section 7, which goes back to being very
expository, deals with the `\'etale' version of the spectral sequence
and the analog of the Quillen-Lichtenbaum conjecture.  Finally, in
section 8 we give the proof of Theorem~\ref{th:main}(a).

\vf

\vfill\eject

\section{Background}

\subsection{Basic setup}

Throughout this paper we will be working in the world of equivariant
homotopy theory over a finite group $G$ (usually with $G=\Zt$).
Unless otherwise indicated, `space' means `equivariant space' and
`map' means `equivariant map'.  If $X$ and $Y$ are spaces, then
$[X,Y]$ denotes the set of equivariant homotopy classes of maps.  When
$H$ is a subgroup of $G$, $[X,Y]^H$ is the set of $H$-equivariant
homotopy classes of $H$-equivariant maps; in particular, $[X,Y]^e$ is
the set of non-equivariant homotopy classes.  The phrase `weak
equivalence' means `equivariant weak equivalence': this refers to a
map $X\ra Y$ such that $X^H \ra Y^H$ is an ordinary weak equivalence
for every subgroup $H\subseteq G$.

\subsection{Connectivity}
\label{se:conn}
Let $V$ be an orthogonal $G$-representation.  Waner \cite[Section
2]{W} introduced the notion of an equivariant space being {\it
$V$-connective\/}, generalizing the non-equivariant notion of
$n$-connectivity.  The key observation is that one can make sense of
the set $[S^{V+k}\Smash G/H_+,X]_*$ not just for $k\geq 0$, but for
$k\geq -|V^H|$ (here, and elsewhere, $|W|$ denotes the real dimension
of the vector space $W$).  If $V_H$ denotes $V$ regarded as an
$H$-representation, then there is a decomposition $V_H=V(H)\oplus
V^H$, where $V(H)$ is the orthogonal complement of the fixed space
$V^H$.  One then considers the chain of equalities
\[ [S^{V+k} \Smash G/H_+,X]_* \cong
   [S^{V_H+k},X]_*^H \cong
   [S^{V(H)+|V^H|+k},X]_*^H 
\]
and observes that the right-hand set makes sense for $k\geq -|V^H|$.
We can therefore take this as a {\it definition\/} for the left-hand
set when $k$ is negative.

A pointed $G$-space $X$ is called \mdfn{$V$-connected} if
$[S^{V+k}\Smash G/H_+,X]_*=0$ for all subgroups $H$ and all $0\geq k\geq
-|V^H|$.  Waner proved that this is equivalent to requiring that $X^H$
is $|V^H|$-connected for all subgroups $H$.  This result eventually
appeared, in expanded form, in \cite[Lemma 1.2]{Lw2}.

The following result of Lewis \cite[Lemma 3.7]{Lw3} will be used
often:

\begin{lemma}
\label{le:lewis}
Suppose $V\supseteq 1$ (the trivial representation), and
let $X$ and $Y$ be pointed $G$-spaces which are both $(V-1)$-connected.
Then a map $X\ra Y$ is a weak equivalence if and only if for every
$k\geq 0$ and every subgroup $H$ it induces an
isomorphism
$[S^{V+k}\Smash G/H_+,X]_{*} \iso [S^{V+k}\Smash G/H_+,Y]_{*}$.
\end{lemma}

\subsection{Eilenberg-MacLane spectra}
\mbox{}\par

\noindent
When $G$ is a finite group, let $\orb(G)$ denote the {\it orbit
category} of $G$---the full subcategory of $G$-spaces whose objects
are the orbits $G/H$.  Recall that a {\it Mackey functor} for $G$ is a
pair of functors $(M^*,M_*)$ from $\orb(G)$ to Abelian groups having
the properties that
\begin{enumerate}[(a)]
\item $M^*$ is contravariant and $M_*$ is covariant;
\item $M^*(G/H)=M_*(G/H)$ for all $H$;
\item For every $t:G/H \ra G/H$ one has $t_*\circ t^*=id$;
\item The double coset formula holds.
\end{enumerate}
We will not write down what the last condition means in general, but
see \cite[XIX.3]{M}.

The importance of Mackey functors is that if $E$ is an equivariant
spectrum and $X$ is any pointed space, then the assignment $G/H
\mapsto [\Sigma^\infty(G/H_{+} \Smash X), E]$ has a natural structure
of a Mackey functor.  In the case $X=S^V$, this Mackey
functor is denoted $\pim_V(E)$ or $\underline{E}^{-V}$.


In the case $G=\Zt$ the orbit category is quite simple,
having the form
\[  \xymatrix{\Zt\ar@(ur,ul)[]_t\ar[r]^i & e }\]
where $it=i$ and $t^2=id$.
It follows that a Mackey functor for $\Zt$ consists of Abelian groups 
$M(\Zt)$ and $M(e)$ together with restriction and transfer maps
\[ \xymatrix{
M(e) \ar@<0.5ex>[r]^{i^*}  & M(\Zt), \ar@<0.5ex>[l]^{i_*} 
& M(\Zt) \ar@<0.5ex>[r]^{t^*}  & M(\Zt)\ar@<0.5ex>[l]^{t_*} 
}
\]
satisfying the following conditions:
\begin{enumerate}[(i)]
\item (Contravariant functoriality) $(t^*)^2=id$ and $t^*i^*=i^*$;
\item (Covariant functoriality) $(t_*)^2=id$ and $i_*t_*=i_*$;
\item $t_*\circ t^*=id$; 
\item (Double Coset formula) $i^* \circ i_* = id+t^*$.
\end{enumerate}
We will specify a Mackey functor for $\Zt$ by specifying the diagram
\[ \xymatrix{ M(\Z/2)  \ar@(ur,ul)[]_{t^*}\ar@<0.5ex>[r]^-{i_*} 
                 & M(e).\ar@<0.5ex>[l]^-{i^*}
}
\]


\begin{example}\mbox{}\par
\label{ex:mackey}
\begin{enumerate}[(a)]
\item The Mackey functor we will be most concerned with is the {\it
constant coefficient} Mackey functor $\Zc$:
\[ \xymatrix{ \Z  \ar@(ur,ul)[]_{id}\ar@<0.5ex>[r]^-{2} 
                 & \Z.\ar@<0.5ex>[l]^-{id}
}
\]
Such a Mackey functor exists over any finite group $G$, and for any
abelian group in place of $\Z$: the restriction maps are all
identities, and the transfer maps $M(G/H) \ra M(G/K)$ are
multiplication by the index $[K:H]$.
\item We define
$\Zc^{op}$ to be
\[ \xymatrix{ \Z  \ar@(ur,ul)[]_{id}\ar@<0.5ex>[r]^-{id} 
                 & \Z.\ar@<0.5ex>[l]^-{2}
}
\]
\item The {\it Burnside ring} Mackey functor $\Am$ is the one for which
$A(G/H)$ is the Burnside ring of $H$.  For $G=\Zt$ this is
\[ \xymatrix{ \Z  \ar@(ur,ul)[]_{id}\ar@<0.5ex>[r]^-{i_*} 
                 & \Z\oplus\Z.\ar@<0.5ex>[l]^-{i^*}
}
\]
where $i^*(a,b)=a+2b$ and $i_*(a)=(0,a)$.  
\end{enumerate}
\end{example}

Every Mackey functor $M$ has an associated $RO(G)$-graded cohomology
theory denoted $V \mapsto H^V(\blank;M)$ (where $V$ runs over all
orthogonal $G$-representations), which is uniquely characterized by
the properties that
\begin{itemize}
\item
$H^n(G/H;M)=
\begin{cases}
    M(G/H) &\text{\ if\ } n=0,\\
    0      &\text{\ otherwise},
\end{cases}
$
\par\noindent
(here $n$ denotes the trivial representation of $G$ on $\R^n$), and
\item
the restriction maps $H^0(G/K;M)\ra H^0(G/H;M)$ induced by
$i\colon G/H\ra G/K$ coincide with the maps $i^*$ in the Mackey functor.
\end{itemize}
The transfer maps
of the Mackey functor will coincide with the transfer maps in this
cohomology theory (or with the pushforward maps in the associated {\it
homology} theory).  Details are in \cite[Chap. IX.5]{M}.

\subsection{Eilenberg-MacLane spaces}
\mbox{}\par

When $M$ is a Mackey functor, the $V$th space in the $\Omega$-spectrum
for $HM$ is called an Eilenberg-MacLane space of type $K(M,V)$.  Such
spaces are $(V-1)$-connected, and have the properties that
$[S^{V+k}\Smash G/H_+,K(M,V)]=0$ for $k>0$ and the Mackey functor $G/H
\mapsto [S^V\Smash G/H_+,K(M,V)]$ is isomorphic to $M$.  See
\cite[Definition 1.4]{Lw3} for this characterization.
  
When $G=\Z/2$ and $V=\R^p \oplus (\R_{-})^q$, we will usually adopt
the motivic notation $K(M,V)=K(M(q),p+q)$.  Likewise $H^V(\blank;M)$
will be written as either $H^{p+q,q}(\blank;M)$ or
$H^{p+q}(\blank;M(q))$, usually the former.

\subsection{The theory $H^{*,*}(X;\Zm)$}
\mbox{}\par

In this section we set down the basic facts about the cohomology
theory $H^*(\blank;\Zc)$ (which we'll sometimes write $H\Zc$).  We
will need to know its coefficient groups $H^{p,q}(\pt;\Zc)$ and
$H^{p,q}(\Zt;\Zc)$, their ring structures, and the transfer and
restriction maps between them.  These things have certainly been
computed many times over the years, although it's hard to find a
precise reference.  The corresponding facts about the theory
$H^*(\blank;\Am)$ can be found in \cite[Thms. 2.1,4.3]{Lw1}, where
they are attributed to Stong.  The necessary information about
$H^*(\blank;\Zc)$ can be deduced from these with a little bit of work,
although it turns out to be much easier to avoid $H^*(\blank;\Am)$
altogether.  The corresponding information about $H^*(\blank;\Zc/2)$
is in \cite[Prop. 6.2]{HK1}, and again one can deduce the integral
analogs with a little bit of work.  An interesting computation of the
{\it positive\/} part of the coefficient ring can be found in
\cite[Thm. 4.1]{LLM}; this gives explicit cycles representing each
element.  In any case, the ultimate conclusions are listed in the
theorem below.  Although the results are not new, we have included
proofs in Appendix B for the reader's convenience.

\begin{thm}\mbox{}\par
\label{th:coeffs}
\begin{enumerate}[(a)]
\item
The abelian group structure of $H^{*,*}(\pt;\Zc)$ is
\[ H^{p,q}(\pt;\Zc)=
\begin{cases}
   \Zt &\text{if $p-q$ is even and $q\geq p >0$;} \\
   \Z &\text{if $p=0$ and $q$ is even;} \\
   \Zt &\text{if $p-q$ is odd and $q+1 < p \leq 0$}\\
   0  &\text{otherwise.}
\end{cases}
\]
These groups are shown in the following picture, where hollow circles
denote $\Z$'s and solid dots represent $\Zt$'s (note that the
$p$-axis is the vertical one):\par
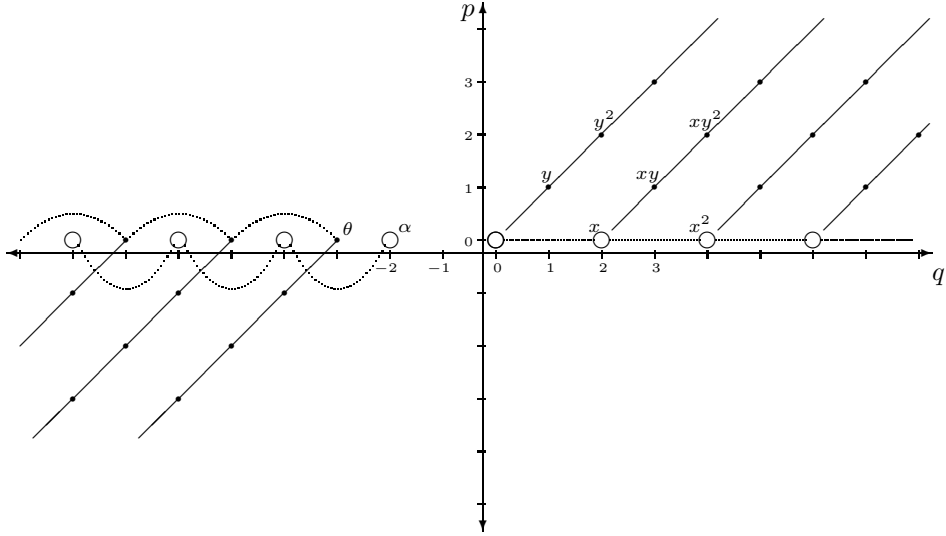
\begin{picture}(330,230)(-150,-110)
\put(-175,5){\vector(1,0){350}}
\put(175,5){\vector(-1,0){350}}
\put(5,-100){\vector(0,1){200}}
\put(5,100){\vector(0,-1){200}}
\multiput(10,10)(40,0){4}{\circle{6}}
\multiput(30,30)(20,20){3}{\circle*{2}}
\multiput(70,30)(20,20){3}{\circle*{2}}
\multiput(110,30)(20,20){3}{\circle*{2}}
\multiput(150,30)(20,20){2}{\circle*{2}}
\multiput(14,14)(40,0){3}{\line(1,1){80}}
\put(134,14){\line(1,1){40}}
%
\multiput(10,10)(-40,0){5}{\circle{6}}
\multiput(-50,10)(-20,-20){4}{\circle*{2}}
\multiput(-90,10)(-20,-20){4}{\circle*{2}}
\multiput(-130,10)(-20,-20){2}{\circle*{2}}
\multiput(-50,10)(-40,0){2}{\line(-1,-1){75}}
\put(-130,10){\line(-1,-1){40}}
\qbezier[30](-68,8)(-50,-25)(-32,8)
\qbezier[30](-108,8)(-90,-25)(-72,8)
\qbezier[30](-148,8)(-130,-25)(-112,8)
%
\qbezier[30](13,10)(30,10)(47,10)
\qbezier[30](53,10)(70,10)(87,10)
\qbezier[30](93,10)(110,10)(127,10)
\qbezier[30](133,10)(150,10)(167,10)
\qbezier[30](-90,10)(-70,30)(-50,10)
\qbezier[30](-130,10)(-110,30)(-90,10)
\qbezier[30](-170,10)(-150,30)(-130,10)
\multiput(-170,3)(20,0){18}{\line(0,1){3}}
\multiput(3,-90)(0,20){10}{\line(1,0){3}}
\put(9,-2){$\scriptscriptstyle{0}$}
\put(29,-2){$\scriptscriptstyle{1}$}
\put(49,-2){$\scriptscriptstyle{2}$}
\put(69,-2){$\scriptscriptstyle{3}$}
\put(-16,-2){$\scriptscriptstyle{-1}$}
\put(-36,-2){$\scriptscriptstyle{-2}$}
\put(-2,8){$\scriptscriptstyle{0}$}
\put(-2,28){$\scriptscriptstyle{1}$}
\put(-2,48){$\scriptscriptstyle{2}$}
\put(-2,68){$\scriptscriptstyle{3}$}
%
\put(45,13){$\scriptstyle x$}
\put(83,13){$\scriptstyle {x^2}$}
\put(27,33){$\scriptstyle y$}
\put(47,53){$\scriptstyle y^2$}
\put(63,33){$\scriptstyle {xy}$}
\put(83,53){$\scriptstyle {xy^2}$}
\put(-27,13){$\scriptstyle \alpha$}
\put(-48,12){$\scriptstyle \theta$}
\put(-3,95){${p}$}
\put(175,-5){${q}$}
\end{picture}
\par\noindent
(An easy way to keep track of the grading is to remember that $y\in
H^{1,1}$ and $x\in H^{0,2}$).
\item The multiplicative structure is completely determined by the
properties that
\begin{enumerate}[(i)]
\item It is commutative;
\item The solid lines in the above diagram represent multiplication by
the class $y\in H^{1,1}$;
\item The dotted lines represent multiplication by $x\in H^{0,2}$ (but
note that only a representative set of dotted lines have been drawn);
\item $x\alpha=2$.
\end{enumerate}
In particular, the subring consisting of $H^{p,q}$ where $p,q\geq 0$
is the polynomial algebra $\Z[x,y]/(2y)$.
%
%
\item The ring $H^{*,*}(\Zt;\Zc)$ is isomorphic to 
$\Z[u,u^{-1}]$, where $u$ has degree $(0,1)$.
\item The Mackey functor $\underline{H\Z}^{0,2n}$ is $\Zm$ when $n\geq
0$ and $\Zm^{op}$ when $n<0$ (see Example~\ref{ex:mackey} for
notation).
\end{enumerate} 
\end{thm}

\begin{remark}
For computations it's often convenient to give every element of
$H^{*,*}(\pt)$ a name in terms of $x$, $y$, and $\theta$.  
For instance, $\alpha$ can be named as $\frac{2}{x}$, and the class
in degree $(-2,-7)$ can be named $\frac{\theta}{xy^2}$.
\end{remark}

\begin{remark}
If $q> 0$ and you look in degree $(1-q,-q)$ and read vertically
upwards, you are seeing the groups $\tH_{sing}^*(\RP^{q-1})$.  
If you look in degree $(q-1,q)$ and read vertically downwards, you are
seeing the groups $\tH_*(\RP^{q-1})$.  The connection is explained in
detail in Appendix B.
\end{remark}

\begin{remark}
It's worth pointing out what aspects of the above picture are similar
to the motivic setting, and which are not.  In the motivic setting one
has that $H^{p,q}_{mot}(\pt;\Zc)=0$ for $q<0$, but notice that this is
{\it not} the case for the above $\Zt$-equivariant theory.  This
difference is tied to the fact that classical algebraic $K$-theory is
connective, whereas $KO$-theory is not.  The Beilinson-Soul\'e
conjecture is that $H^{p,q}_{mot}(\pt;\Zc)=0$ when $q>0$ and $p<0$,
which is clearly satisfied in our $\Zt$-world.  The
non-zero motivic cohomology groups of a point should correspond to the
groups lying in the first quadrant of the above diagram, with the same
vanishing line.  The motivic groups lying {\it along} this line are the
Milnor $K$-theory groups.
\end{remark}

The above result tells us everything about the Eilenberg-MacLane
spaces $K(\Z(n),2n)$.  Here are a couple of the main points:

\begin{cor}\mbox{}\par
\begin{enumerate}[(a)]
\item Non-equivariantly, $K(\Z(n),2n)$ is a $K(\Z,2n)$.
\item The induced action of $\Zt$ on $\pi_{2n}(K(\Z(n),2n))=\Z$ is
multiplication by $(-1)^n$.
\item The fixed set $K(\Z(n),2n)^{\Z_2}$ has the homotopy type of
either
\[ K(\Z,2n)\times K(\Zt,2n-2)\times K(\Zt,2n-4)\times \cdots \times
K(\Zt,n) \qquad (n \text{\ even,}) \]
or
\[ K(\Zt,2n-1)\times K(\Zt,2n-3)\times \cdots \times K(\Zt,n) \qquad
(n \text{\ odd}).\]
\end{enumerate}
\end{cor}

\begin{proof}
A theorem of \cite{dS} identifies $K(\Z(n),2n)$ with $AG(S^{2n,n})$,
the free abelian group generated by $S^{2n,n}$.  So both $K(\Z(n),2n)$
and its fixed set are topological abelian groups, hence products of
Eilenberg-MacLane spaces.  The homotopy groups can be read off of
$H^{*,*}(\pt;\Zc)$ and $H^{*,*}(\Zt;\Zc)$.  This proves (a) and (c).

Note that 
$[S^{2n,n},K(\Z(n),2n)]^e\iso [S^{2n,n}\Smash
\Zt_+,K(\Z(n),2n)]_* \iso H^{0,0}(\Zt)$, 
and we know the group action on the latter is trivial (because we know
the Mackey functor $\underline{H}^{0,0}$).  The group
$\pi_{2n}K(\Z(n),2n)$ may be written $[S^{2n,0},K(\Z(n),2n)]^e$,
and this differs from the above in the replacement of $S^{2n,n}$ by
$S^{2n,0}$.  On the former sphere, the automorphism coming
from the $\Zt$ action has degree $(-1)^n$ (complex conjugation on
$\C^n$ reflects $n$ real coordinates).  This proves (b).
\end{proof}

\begin{remark}
\label{re:EMhfixed}
The homotopy fixed set
$K(\Z(n),2n)^{h\Z_2}$ will also be a topological abelian group, and
hence a generalized Eilenberg-MacLane space.  The spectral sequence
for computing homotopy groups of a homotopy limit collapses, and shows
that $K(\Z(n),2n)^{h\Zt}$ is either
\[ K(\Z,2n)\times K(\Zt,2n-2)\times K(\Zt,2n-4)\times \cdots \times
K(\Zt,0) \qquad (n \text{\ even}), \]
or
\[ K(\Zt,2n-1)\times K(\Zt,2n-3)\times \cdots \times K(\Zt,1) \qquad
(n \text{\ odd}).\]
So the actual fixed set is a truncation of the homotopy fixed
set.  This observation reappears in section~\ref{se:etale}.
\end{remark}

\vf

\section{Equivariant Postnikov-section functors}

In this section we define two types of equivariant Postnikov section
functors, denoted $\PK_V$ and $\PP_V$, and list their basic properties.  

\medskip

To begin this section we work in the context of an arbitrary finite
group $G$.  The category $\Gtop$ denotes the category of $G$-spaces
which are compactly-generated and weak Hausdorff, with equivariant
maps.  We will eventually specialize to the case $G=\Z/2$, but for the
present it is just as easy to work in greater generality.  There is a
model category structure on $\Gtop$ analagous to the usual one on
$\Top$.  

\subsection{Generalities}
Recall that a space $A$ is said to be \mdfn{small with respect to
closed inclusions} if it has the property that for any sequence
of closed inclusions
\[ Z_0 \inc Z_1 \inc Z_2 \inc \cdots \]
the canonical map $\colim_i \Gtop(A,Z_i) \ra \Gtop(A,\colim Z_i)$ is
an isomorphism.  Every compact Hausdorff space is small in this sense.
Let $CA$ denote the cone on $A$, and recall that $[X,Y]$ denotes {\it
unpointed\/} homotopy classes of maps.

Let $\cA$ be a set of well-pointed spaces, all of which are compact
Hausdorff.  The pointedness can be ignored for the moment, but will
be needed later.  We will say that a space $Z$ is \mdfn{$\cA$-null} if
it has the property that the maps $[*,Z] \ra [\Sigma^n A,Z]$ (induced
by $\Sigma^n A \ra *$) are isomorphisms, for all $n\geq 0$ and all
$A\in \cA$.  This is equivalent to saying that every map $\Sigma^n
A\ra Z$ extends over the cone.

For a space $X$ one can construct a new space $P_\cA(X)$ with the
following properties:
\begin{enumerate}[(1)]
\item There is a natural map $X \ra P_{\cA}X$;
\item $P_{\cA}X$ is $\cA$-null;
\item If $Z$ is an $\cA$-null space and $X\ra Z$ is a map, then
there is a lifting
\[ \xymatrix{ X \ar[r]\ar[d] & Z \\
              P_{\cA}X \ar@{.>}[ur]
}\]
and this lifting is unique up to homotopy.
\end{enumerate}
The functors $P_\cA(X)$ are called \dfn{nullification functors} in
\cite{F}.  They are examples of Bousfield localization functors, 
for which an excellent reference is \cite[Chapters 3,4]{H}.
In our context we construct them as follows:
For any space $Y$, let $F_{\cA}(Y)$ be defined by the 
pushout square
\[ \xymatrix{ 
         \coprod_{\sigma} \uSigma^n A \ar[r]\ar@{ >->}[d] 
                       & Y \ar@{ >.>}[d] \\
         \coprod_{\sigma} C(\uSigma^n A) \ar@{.>}[r] & F_{\cA} Y, 
}
\]
where $\sigma$ runs over all maps $\uSigma^n A \ra Y$ (for all $A\in
\cA$, and $n\geq 0$).  One then considers the sequence of closed
inclusions
\[ X \inc F_{\cA}X \inc F_{\cA}F_{\cA}X \inc F_{\cA}F_{\cA}F_{\cA}
X \inc \cdots \]
and $P_{\cA}(X)$ is defined to be the colimit.  It is routine to
check, using basic obstruction theory, that this construction has the
required properties.

\begin{remark}
We will often use the observation that if $\cA \subseteq \cA'$ then
one has a canonical map $P_{\cA}X \ra P_{\cA'}X$.
\end{remark}

In terms of the localization of model categories (cf. \cite{H}), we
are localizing $\Gtop$ at the set of maps $\{A \ra *| A\in \cA\}$ and
$P_{\cA}$ is the localization functor.  This uses the fact that the
objects in $\cA$ are well-pointed.  As a consequence, the functors
$P_{\cA}X$ have the standard properties one would expect from a
localization functor.  We omit the proof of the following result: such
properties can be found in \cite{H} in complete generality, or in
\cite{F} for the category of spaces.  In this case they are also easy
to prove directly by standard arguments.


\begin{prop} \mbox{}\par
\label{pr:locohoco}
\begin{enumerate}[(a)]
\item
Let $X\ra Y$ and $X\ra Z$ be maps, where $X\ra Y$ is a cofibration.
If $Z$ is $\cA$-null and
$P_{\cA}X\ra P_{\cA}Y$ is a weak equivalence, then there is a lift
\[ \xymatrix{ X \ar[d]\ar[r] & Z\\
              Y \ar@{.>}[ur]
}\]
and this lift is unique up to homotopy.
\item
Let $X:\cC\ra \Gtop$ be a diagram.  Then the natural map
\[ P_{\cA}(\hocolim_\alpha X_\alpha) \lra P_{\cA}(\hocolim P_{\cA}X_\alpha) \]
is a weak equivalence.
\item 
If $X\ra Y\ra Z$ is a homotopy cofiber sequence and $P_{\cA}X$ is
contractible, then $P_{\cA}Y\ra P_{\cA}Z$ is a weak equivalence.
\end{enumerate}
\end{prop}

\subsection{The functors $\PK$ and $\PP$}
\label{se:postnikov}

The most basic examples of nullification functors are the ordinary
Postnikov section functors (when $G=\{e\}$): given a non-equivariant
space $X$ one forms $P_n X$ by killing off all homotopy groups above
dimension $n$.  In the language of the previous section $P_n
X=P_{\cA_n} X$, where $\cA_n$ is the set $\{S^{n+1},S^{n+2},\cdots
\}$.  In fact we would get a homotopy equivalent space by just taking
$\cA_n$ to be $\{S^{n+1}\}$, but we have arranged things so that
$\cA_{n+1} \subseteq \cA_{n}$ because our formalism then gives natural
maps $P_{\cA_{n+1}}X \ra P_{\cA_n}X$.

When one wants to introduce Postnikov section functors for $G$-spaces
several possibilities present themselves.  One thing to note is that
when we kill off all maps from a space $A$, we would also like to be
killing off maps from all spaces $A\Smash Z$.  In the nonequivariant
setting this is automatic, because $Z$ can be built from spheres
and therefore $A\Smash Z$ is built from suspensions of $A$.  
In the equivariant setting we have to explicitly build
this into the theory, by making sure that whenever we kill off a space
$A$ we also kill off $A\Smash G/H_+$ for all subgroups $H$.

If $V$ is a representation of $G$ then one Postnikov functor we can
consider is $P_\cA(X)$ where $\cA$ is the set
$\{S^{V+n}\Smash G/H_+ | n> 0, H\leq G\}$.
One can get the same result by
doing the following, and for functorial reasons it is somewhat better:
Let 
\[ \tilde{\cA}_{V}= \{ S^W \Smash G/H_+ \mid W \supseteq V+1, H\leq G \}, \]
and define $\PK_V X = P_{\tilde{\cA}_V}(X)$.  This definition
guarantees that if $V \subseteq U$ then there are natural maps
$\PK_U X \ra \PK_V X$.

On the other hand we can also do the following. Let
\[ \cA_V= \{S^W \Smash G/H_+ \mid W \supset V, H\leq G \}, \]
and define $\PP_V X=P_{\cA_V}(X)$.  Again, whenever $V\subseteq U$
there are natural maps $\PP_U X \ra \PP_V X$.  Moreover, since
$\tilde{\cA}_V \subset \cA_V$ one has maps $\PK_V X \ra \PP_V X$.

\begin{remark}
The difference between $\PK_V$ and $\PP_V$ shows up in the following
way.  Non-equivariantly, there are no non-trivial maps $S^n \ra S^k$
when $k>n$.  As a result, when one forms the Postnikov section $P_n X$
one doesn't change the homotopy groups in low dimensions: $[S^m,X] \ra
[S^m,P_n X]$ is an isomorphism for $m\leq n$.  In the equivariant
theory, however, there can be many non-trivial maps $S^V \ra S^W$ for
$V\subset W$; so when forming a Postnikov section by killing off
maps from large spheres, one may actually be creating {\it new\/} maps
from smaller spheres.

What is true, however, is that if $V+1\subseteq W$---that is, if $W$
contains $V$ plus at least one copy of the trivial
representation---then all equivariant maps $S^V \ra S^W$ are null.
This leads one to the functors $\PK_V$ defined above, which are
designed so that they don't change homotopy classes of maps from $S^V$
and smaller spheres.  The general rule is that the functors $\PK_V$
are better behaved that $\PP_V$: they are easier to compute, and their
properties (outlined below) closely resemble those of non-equivariant
Postnikov sections.  These are the same as the Postnikov functors in
\cite[II.1]{M}.
\end{remark}
 
\begin{prop}[Properties of $\PK$] 
\label{pr:pkprop}
If $X$ is a pointed $G$-space $X$ and $V$ is a
$G$-representation, the following are true:
\begin{enumerate}[(a)]
\item The map $X \ra \PK_V X$ induces an isomorphism of the sets
$[S^{k,0}\Smash G/H_+, \blank]_*$ for $0 \leq k \leq \dim V^H$,
and an epimorphism for $k=\dim V^H+1$.
\item If $W$ is a $G$-representation for which $\dim W^H \leq \dim
V^H$ for all subgroups $H\subseteq G$, then $[S^W,X]_* \ra
[S^W,\PK_V X]_*$ is an isomorphism.
\item The homotopy fiber of $\PK_{V+1} X \ra \PK_V X$ is an 
Eilenberg-MacLane space of type $K(\upi_{V+1} X, V)$.
\item The homotopy limit of the sequence 
\[ \cdots \ra \PK_{V+2} X \ra \PK_{V+1}X \ra \PK_V X \]
is weakly equivalent to $X$.
\item If $V$ contains the regular representation of $G$, then the
Postnikov section $\PK_V(S^V)$ is an Eilenberg-MacLane space of type
$K(\Am,V)$ where $\Am$ is the Burnside-ring Mackey functor.
\end{enumerate}
\end{prop}

\begin{proof}
The proof is completely standard, so we will only give a brief sketch.
Suppose that $X$ is a space, $S^W \Smash G/J_+ \ra X$ is a map, and we
attach the cone on this map to construct a new space $X_1$.  The for
any subgroup $H\leq G$, $X_1^H$ is obtained from $X^H$ by attaching a
cone on the map $(S^W \Smash G/J_+)^H \ra X^H$.  The domain of this
map is a wedge of spheres of dimension $|W^H|$, and so $X^H \ra
X_1^H$ is $|W^H|$-connected.  From these considerations (a) is
immediate.

Part (b) follows from (a) and the fact that $S^W$ has an equivariant
CW-structure made up of cells $S^{k,0}\Smash G/H_+$ where $k\leq \dim
W^H$.  Part (d) is also immediate from (a): the map $X\ra \PK_V X$
becomes highly connected on all fixed sets as $V$ gets large.

For part (c), note first that for an arbitrary space $X$ the object
$\underline{\pi}_V X$ may not be a Mackey functor---it is instead a {\it
$V$-Mackey functor} as defined in \cite[1.2]{Lw3}.  A characterization
of Eilenberg-MacLane spaces is given in \cite[1.4]{Lw3}, and it is
easy to use parts (a) and (b) to check that the homotopy fiber we're
looking at has the properties listed there.

Part (d) is an immediate consequence of (c) and the well-known
isomorphism of Mackey functors $\upi_V(S^V)\iso \Am$ (which follows
from \cite[IX.1.4,XVII.2]{M}).
\end{proof}

Suppose now that $\cE$ is an equivariant spectrum, and let $E$ denote
the $0$th space of the corresponding $\Omega$-spectrum.  By applying
the functors $\PK_n$ we obtain a tower
\[ \xymatrix{
  \cdots \ar[r] & \PK_2 E \ar[r] & \PK_1 E \ar[r] & \PK_0 E \ar[r] & {*}\\
& K(\upi_2 E,2) \ar[u] &K(\upi_1 E,1) \ar[u] & K(\upi_0 E,0) \ar[u] \\
}
\]
The homotopy spectral sequence for maps from $X$ into this tower gives
the classical equivariant Atiyah-Hirzebruch spectral sequence
$H^p(X;\underline{\cE}^q) \Rightarrow \cE^{p+q,0}(X)$ with a suitable
truncation.  Here $\underline{\cE}^q$ denotes the Mackey functor $G/H
\mapsto \cE^q(G/H)$.  To get the full spectral sequence one can look
at the Postnikov towers for each $\cE_V$ (the $V$th space in the
$\Omega$-spectrum for $\cE$) and note that the resulting spectral
sequences can be pasted together.

Unfortunately, for $\cE=\KR$ this spectral sequence is not the one
we're looking for.  The above spectral sequence collapses when $X=*$
and gives no information, whereas the spectral sequence we're looking
for is very non-trivial when $X=*$.  In the context of algebraic
$K$-theory, the analog of the above spectral sequence is the
{\it Brown-Gersten spectral sequence\/} of \cite{BG}.  
 
\medskip

The functors $\PP_{V}$ don't have the property that the homotopy fiber
of $\PP_{V+1}X \ra \PP_V X$ is necessarily an equivariant
Eilenberg-MacLane space.  For the record, here are the basic
properties of $\PP$.  The proofs are the same as for
Proposition~\ref{pr:pkprop}.  

\begin{prop}[Properties of $\PP$] 
\label{pr:ppprop}
For any pointed $G$-space $X$ and any
$G$-representation $V$, the following are true:
\begin{enumerate}[(a)]
\item The map $X \ra \PP_V X$ induces an isomorphism of the sets
$[S^{k,0}\Smash G/H_+, \blank]_*$ for $0 \leq k < \dim V^H$,
and an epimorphism for $k=\dim V^H$.
\item If $W$ is a $G$-representation for which $\dim W^H < \dim
V^H$ for all subgroups $H\subseteq G$, then $[S^W,X]_* \ra
[S^W,\PK_V(X)]_*$ is an isomorphism.
%
%
\end{enumerate}
\end{prop}

The main reason we care about the functors
$\PP_V$ is the following result:

\begin{thm}
\label{th:P(S)}
When $G=\Zt$ and $V$ contains the trivial representation,
the space $\PP_V(S^V)$ has the equivariant weak homotopy type of the
Eilenberg-MacLane space $K(\Zm,V)$.
\end{thm}

The proof of this result is somewhat involved, and will be postponed
until section~\ref{se:symmprod}.  However, we can give some intuition
for why it's true.  If $V\supseteq \C$ one knows that
$[S^V,S^V]_*=\Z\oplus \Z$ (the Burnside ring of $\Zt$), and $[S^V,S^V]_*
\ra [S^V,\PK_V(S^V)]_*$ is an isomorphism by
Proposition~\ref{pr:pkprop}(b).  If one chooses the generators of
$[S^V,S^V]_*$ appropriately, their difference factors through a `Hopf
map' $S^{V+\R_{-}} \ra S^V$, where $\R_{-}$ is the sign representation of
$\Zt$.  Since $\PP_V(S^{V+\R_{-}})\he *$, this difference becomes null in
$\PP_V(S^V)$ (note that $\PK_V(S^{V+\R_{-}})$ is not contractible).  So
the two copies of $\Z$ in $[S^V,\PK_V(S^V)]_*$ become identified in
$[S^V,\PP_V(S^V)]_*$, and this is ultimately why $\PP_V(S^V)$ has the
homotopy type of $K(\Zm,V)$ rather than $K(\Am,V)$.

We'd like to justify the claim that the difference of generators
factors through a Hopf map.  First look at the case $V=\C$, so that
$S^V=S^{2,1}$.  Consider the degree map $[S^{2,1},S^{2,1}]_* \ra
\Z\oplus \Z$ which sends a map $f$ to the pair $(\deg f, \deg
f^{\Zt})$ (the degree, and the degree of the map restricted to the
fixed set).  This is injective, and the image consists of pairs
$(n,m)$ for which $n-m\equiv 0 \mod 2$.  For the generators of
$[S^{2,1},S^{2,1}]_*$ we'll take the identity map, which has degree
$(1,1)$, and the complex conjugation map, which has degree $(-1,1)$.

Now consider the Hopf projection $\C^2-\{0\} \ra P_{\C}^1$, which we
may write as $p\colon S^{3,2}\ra S^{2,1}$.  Smashing the inclusion
$S^{0,0} \inc S^{1,1}$ with $S^{2,1}$ gives a map $j\colon S^{2,1}
\inc S^{3,2}$, and the degree of the composition $pj$ is
$(0,2)$.  So this composite is homotopic to the sum of our two
chosen generators.  When $V=\C\oplus W$, one takes this same argument
and smashes everything in sight with $S^W$.

Another perspective on Theorem~\ref{th:P(S)} is given in
Section~\ref{se:symmprod}, where it is tied to the geometry of the
infinite symmetric product construction.
  
\vf


\section{The Postnikov tower for $\Z\times BU$}

In this section we will consider the objects $\PP_{n\C}(\Z\times BU)$.
In motivic indexing $\PP_{n\C}$ would be written $\PP_{(2n,n)}$, and
we will abbreviate this as just $\PP_{2n}$.  Our goal is the following

\begin{thm} 
\label{th:unstablefiberseq}
Let $\beta:S^{2,1}\ra\Z\times BU$ be a map representing
the Bott element in $\widetilde{KR}^{0,0}(S^{2,1})$, and let 
$\beta^n:S^{2n,n}\ra \Z\times BU$ denote its $n$th power.  
Then
\[ \PP_{2n}(S^{2n,n}) \llra{\beta^n} \PP_{2n}(\Z\times BU) \lra
\PP_{2n-2}(\Z\times BU) 
\]
is a homotopy fiber sequence.   
\end{thm}

\begin{cor}  
\label{co:ptower}
There is a tower of homotopy fiber sequences
\[
   \xymatrix{   
                \quad \cdots\quad \ar[r]
                &\PP_{4}(\KK) \ar[r]  
                &\PP_{2}(\KK) \ar[r] 
                &\PP_{0}(\KK) \ar[r]
                &{*}\\            
            &K(\Z(2),4)\ar[u]
            &K(\Z(1),2)\ar[u]
            &\Z \ar[u]
}
\]
and the homotopy limit of the tower is $\KK$.
\end{cor}

\smallskip

As explained in the introduction, one can prove the corollary using
only the functors $\PK$, and this is easier in the end.
Writing $\PK_{2n}$ for $\PK_{n\C}$, we have the following:

\begin{prop}
The homotopy fiber of $\PK_{2n}(\KK) \ra \PK_{2n-2}(\KK)$ is an
Eilenberg-MacLane space of type $K(\Z(n),2n)$.
\end{prop}

\begin{proof}
Let $\F_n$ denote the homotopy fiber.  From
Proposition~\ref{pr:pkprop}(a) it follows immediately that the fixed
set of $\F_n$ is $(n-1)$-connected.  Non-equivariantly $\PK_{2n}$ has
the homotopy type of the ordinary Postnikov section functor, and so
$\F_n$ is $(2n-1)$-connected as a non-equivariant space (here we are
using that $\Z\times BU$ has no odd homotopy groups).  So in
equivariant language $\F_n$ is $(n\C-1)$-connected
(cf. section~\ref{se:conn}).  By construction we have that
$[S^{2n+k,n}\Smash \Zt_+,\F_n]_*=0=[S^{2n+k,n},\F_n]_*$ for all $k>
0$, and the long exact homotopy sequence shows that the Mackey functor
$\underline{\pi}_{2n,n}(\F_n)$ is isomorphic to
$\underline{\pi}_{2n,n}(\KK)\iso \underline{KR}^{2n,n}(\pt) \iso \Zm$.
Lewis's characterization of Eilenberg-MacLane spaces \cite[Def. 1.4
]{Lw3} now shows that $\F_n$ is a $K(\Z(n),2n)$.
\end{proof}

\smallskip

\noindent
For the remainder of the section we let $F_n$ denote the homotopy
fiber  
\[ F_n \lra \PP_{2n}(\Z\times BU) \lra \PP_{{2n-2}}(\Z\times BU).\]
The map $\beta^n\colon \PP_{2n}(S^{2n,n}) \ra \PP_{2n}(\KK)$ becomes
null when we pass to the space $\PP_{2n-2}(\KK)$, and therefore it
lifts to $F_n$ (and the lifting is unique up to homotopy).  Our task
is to show that this lifting is a weak equivalence.  Both the domain
and codomain are highly connected, and so we can use
Lemma~\ref{le:lewis}.  Here is a restatement for the special case
where $G=\Zt$ and $V=\C^n$:

\begin{lemma}
\label{le:waner}
Let $X$ and $Y$ be pointed $\Zt$-spaces with the properties that
\begin{enumerate}[(i)]
\item $[S^{k,0},X]_*=[S^{k,0},Y]_*=0$ for $0\leq k<n$, and
\item $[\Zt_+\Smash S^{k,0},X]_*=[\Zt_+\Smash S^{k,0},Y]_*=0$ for $0\leq k
< 2n$. 
\end{enumerate}
Then a map $X\ra Y$ is an equivariant weak equivalence if and only if
it induces isomorphisms
\[ 
   [S^{2n+k,n},X]_* \mapisot [S^{2n+k,n},Y]_* \ \text{and}\ \  
   [\Zt_+\Smash S^{2n+k,n},X]_* \mapisot [\Zt_+\Smash S^{2n+k,n},Y]_* 
\]
for every $k\geq 0$.
\end{lemma}

We know by Theorem~\ref{th:P(S)} that $\PP_{2n}(S^{2n,n})$ is a
$K(\Z(n),2n)$-space and therefore we know it's homotopy groups---these
are precisely the groups $H^{p,q}(\pt;\Zc)$ and $H^{p,q}(\Zt;\Zc)$.
So it's easy to see that $K(\Z(n),2n)$ satisfies the conditions in the
above lemma.  The general strategy at this point would be to
\begin{enumerate}[(a)]
\item Show that $F_n$ also satisfies the conditions of the lemma;
\item Observe that $[S^{2n+k,n},F_n]_*=0=[\Zt_+\Smash S^{2n+k,n},F_n]_*$ 
      for $k>0$, for trivial reasons;
\item Show that the map $\PP_{2n}(S^{2n,n})\ra F_n$ induces isomorphisms on
      $[S^{2n,n},\blank]_*$ and $[\Zt_+\Smash S^{2n,n},\blank]_*$;
\item Use Lemma~\ref{le:waner} to deduce that 
      $\PP_{2n}(S^{2n,n})\ra F_n$ is a weak equivalence.
\end{enumerate}

In fact this approach can be streamlined a bit by using the functors
$\PK$ as a crutch.  

\begin{lemma}\mbox{}\par
\label{le:tildePBU}
\begin{enumerate}[(a)]
\item Let $X$ be a pointed $\Zt$-space with the property that the 
      forgetful map $[S^{2n,n},X]_*\ra [S^{2n},X]^e_*$ is injective.  
      Then the natural map $\PK_{2n} X\ra \PP_{2n} X$ is a weak 
      equivalence.
\item $\PK_{2n}(\Z\times BU) \ra \PP_{2n}(\Z\times BU)$ is a
      weak equivalence.
\end{enumerate}
\end{lemma}

\begin{proof}
For (a) we only have to show that $[S^{2n+p,n+p},\PK_{2n}X]=0$ for
all $p>0$; that is, we must show that $\PK_{2n} X$ is null with
respect to $\cA_{(2n,n)}$, not just $\tilde \cA_{(2n,n)}$ (see
section~\ref{se:postnikov}).  Consider the basic Puppe sequence 
$\Zt_+ \ra S^{0,0} \ra S^{1,1} \ra \Zt_+ \Smash S^{1,0}$.
Smashing with $S^{2n,n}$ yields
\[ \Zt_+\Smash S^{2n,n} \ra S^{2n,n} \ra S^{2n+1,n+1} \ra \Zt_+ \Smash
   S^{2n+1,n}.
\]
Mapping this sequence into $\PK_{2n}X$ gives the top edge of the
following diagram:
\[ \xymatrix{
&[S^{2n+1,n+1},\PK_{2n}X] \ar[d]
           & [\Zt_+\Smash S^{2n+1,n},\PK_{2n} X]=0\ar[l]\\
       [\Zt_+\Smash S^{2n,n},\PK_{2n}X] 
           & [S^{2n,n},\PK_{2n}X]\ar[l] \\
       [\Zt_+\Smash S^{2n,n},X] \ar[u]^\cong
           & [S^{2n,n},X]. \ar[u]^\cong\ar[l]
}
\]
The right-most group in the top row is zero just because of the
definition of $\PK_{2n}$, and Proposition~\ref{pr:pkprop}(b) implies
that the labelled vertical maps are isomorphisms.  The map in the
bottom row may be identified with the forgetful map
\[ [S^{2n},X]^e \la [S^{2n,n},X], \]
and we have assumed that this is injective.  It's  now clear that
$[S^{2n+1,n+1},\PK_{2n}X]$ must be zero.

Smashing the above Puppe sequence with $S^{2n+p,n+p}$ gives
\[ S^{2n+p,n+p} \ra S^{2n+p+1,n+p+1} \ra \Zt_+\Smash S^{2n+p+1,n+p}.
\]  
By induction we know that $\PK_{2n}X$ is null with respect to the
first and third space, so it is also null with respect to the second.
This finishes (a).

\smallskip

Proving (b) is of course just a matter of checking that $\Z\times BU$ has the
property specified in (a).  So we must check that the forgetful map
\[ \Z={\widetilde{KR}}\phantom{|}^{0,0}(S^{2n,n}) \ra \tilde{K}^0(S^{2n})=\Z \]
is injective.  But the map is easily seen to be an isomorphism, as 
$\beta^n$ is an explicit generator for both the domain and target.
\end{proof}

\begin{proof}[Proof of Theorem~\ref{th:unstablefiberseq}]
We must show that $j\colon \PP_{2n}(S^{2n,n}) \ra F_n$ is a weak
equivalence.  The equivalences $\PK_{2n}(\KK) \ra \PP_{2n}(\KK)$
induce equivalences $\F_n \ra F_n$, and we already know $\F_n\he
K(\Z(n),2n)$.  Lemma~\ref{le:waner} now implies that $j$ must be a weak
equivalence (one uses the fact that $\beta^n$ is a generator for
$KR^{2n,n}(\pt)$).
\end{proof}

\begin{proof}[Proof of Corollary~\ref{co:ptower}]
This is just a a restatement, together with the fact that the holim of
the tower is $\Z\times BU$.  The latter follows from
Proposition~\ref{pr:ppprop}(a,b).
\end{proof}

\vf


\section{Properties of the spectral sequence}

If $X$ is a $\Zt$-space, then the associated homotopy spectral
sequence for the tower of 
Corollary~\ref{co:ptower} has the form
\[ H^{p,-\frac{q}{2}}(X;\Zc)\Rightarrow 
                     [S^{-p-q,0}\Smash X_+,\Z\times BU]_*,
\]
being confined to the quadrant $p,q\leq 0$.  
This is an unstable version of the spectral sequence we're looking
for. Producing the stable version is not difficult, as one can replace
$X$ by various suspensions $S^{a,b}\Smash X$ and use the periodicity
of $\Z\times BU$ to get a `family' of spectral sequences which patch
together.  We'll take another approach to this in the next section,
and for now be content with analyzing the unstable case.

\subsection{Adams operations}\mbox{}\par

\smallskip

There is a map of $\Zt$-spaces $\psi^k:\KKr \ra \KKr$ inducing the
operation $\psi^k$ on $KR^0(X)$, constructed out of the $\lambda^i$
maps in the usual way.  The functoriality of the constructions
$\PP_{2n}$ shows that $\psi^k$ induces a self-map of the Postnikov
tower for $\KKr$, and therefore we get an action of the Adams
operations on the spectral sequence.  We must identify the action on
the $E_2$-term:

\begin{prop}
\label{pr:adams}
The induced map $\psi^k:F_n \ra F_n$ coincides with the multiplication
by $k^n$ map $K(\Z(n),2n)\ra K(\Z(n),2n)$.
\end{prop}

\begin{proof}
If $\beta^n:S^{2n,n}\ra \Z\times BU$ is the $n$th power of the Bott
element, then we know the following diagram commutes:
\[ \xymatrix{
      S^{2n,n} \ar[r]\ar[d]_{\cdot k^n} & \KKr \ar[d]^{\psi^k} \\
      S^{2n,n} \ar[r] & \KKr.
}\]
This is just because we can compute
\[ \psi^k(\beta^n)=(\psi^k\beta)^n=(k\beta)^n=k^n\beta^n.\]
Applying $\PP_{2n}$ to the above diagram gives
\[ \xymatrix{
      S^{2n,n}\ar[d]_{k^n}\ar[r] &\PP_{2n}(S^{2n,n})
           \ar[r]\ar[d]^{\PP_{2n}(k^n)} & \PP_{2n}(\KKr)
           \ar[d]^{\PP_{2n}(\psi^k)} \\ S^{2n,n}\ar[r] &\PP_{2n}(S^{2n,n})
           \ar[r] & \PP_{2n}(\KKr).  
}
\] 
We have previously identified the map $F_n \ra \PP_{2n}(\KKr)$ with
$\PP_{2n}(S^{2n,n})\ra \PP_{2n}(\KKr)$, and so the argument may be
completed by proving the following lemma.
\end{proof}

\begin{lemma}
Let $k\in \Z$ and let $k:S^{2n,n}\ra S^{2n,n}$ denote the map obtained
by adding the identity to itself $k$ times in the group
$[S^{2n,n},S^{2n,n}]_*$ (using the fact that $S^{2n,n}$ is a
suspension).  Then the localized map
\[ \PP_{2n}(k):\PP_{2n}(S^{2n,n})\ra \PP_{2n}(S^{2n,n}) \]
may be identified with the map $K(\Z(n),2n) \ra K(\Z(n),2n)$
representing multiplication by $k$. 
\end{lemma}

\begin{proof}
There are several ways one could do this.
Write $S$ for $S^{2n,n}$ and $P$ for $\PP_{2n}(S)$.  We of course have the
diagram
\[ \xymatrix{ S \ar[r]\ar[d]_{k} & P\ar[d]^{P(k)}\\
              S \ar[r] & P.
}\]
Using Proposition~\ref{pr:locohoco}(a) it's easy to see that $[P,P]_*\ra
[S,P]_*$ is an isomorphism, and the arguments in
Section~\ref{se:symmprod} show that $[S,P]\iso \Z$ is generated by 
the localization map $S\ra P$.  This proves it.
\end{proof}

\subsection{The rational tower}\mbox{}\par

Grassmannians have nice Schubert cell decompositions, which make it
easy to compute $H^{*,*}(\blank)$.  One of course finds that
$H^{*,*}(BU)=H^{*,*}(\pt)[c_1,c_2,\ldots]$ where $c_i$ has degree
$(2i,i)$.  If we regard $c_n$ as a map $BU \ra K(\Z(n),2n)$, then applying
$P_{2n}$ gives $P_{2n}(BU) \ra P_{2n}K(\Z(n),2n) = K(\Z(n),2n)$ (the
Eilenberg-MacLane space is already $\cA_{(2n,n)}$-null).  We claim the
composite
\[ K(\Z(n),2n) =P_{2n}(S^{2n,n}) \llra{P(\beta^n)} 
P_{2n}(\Z\times BU) \llra{c_n} K(\Z(n),2n) 
\]
is multiplication by $(n-1)!$.  As in the last section, the argument
comes down to knowing that $c_n(\beta^n)$ is $(n-1)!$ times the
generator of $\tH^{2n,n}(S^{2n,n})$.  This can be deduced via comparison
maps to the nonequivariant groups, where the result is well-known (due
to Bott, originally).

So we see that the inclusions of homotopy fibers $K(\Z(n),2n) \ra
P_{2n}(\Z\times BU)$ are split rationally; hence the spectral sequence
collapses rationally.  If $\oplus_n H^{2n,n}(X)\tens \Q$ is
finite-dimensional, then $KR^0(X)\tens \Q$ decomposes into eigenspaces
of the Adams operations.

\subsection{Convergence}\mbox{}\par

The homotopy spectral sequence for a bounded below tower is
automatically conditionally convergent \cite[Def. 5.10]{Bd}.
So if $RE_\infty=0$ it converges strongly, by \cite[7.4]{Bd} and the
Milnor exact sequence.

\subsection{Multiplicativity}\mbox{}\par
\label{se:mult}

The proof of multiplicativity follows the same lines as the
nonequivariant case, which is written up in detail in \cite{D2}.
We will only give an outline.

One starts by letting $W_{n}$ be the homotopy fiber of $\KKr \ra
P_{2n-2}(\KKr)$, and these come with natural maps $W_{n+1} \ra W_n$.
The square
\[ \xymatrix{ \KKr \ar[r] \ar[d] & P_{2n-2}(\KKr) \ar[d] \\
              P_{2n}(\KKr) \ar[r] & P_{2n-2}(\KKr)
}
\]
gives us a map $W_n \ra F_n$, where $F_n$ is the homotopy fiber of the
bottom map (which we know is a $K(\Z(n),2n)$). 
Routine nonsense shows that $W_{n+1} \ra W_n \ra F_n$ is a homotopy
fiber sequence, and so we have a tower 
\[
   \xymatrix{
&
K(\Z(3),6) &
K(\Z(2),4) &
K(\Z(1),2) & 
K(\Z(0),0) \\
                \quad \cdots\quad  \ar[r]
                &W_3 \ar[u]\ar[r]
                &W_2 \ar[u]\ar[r]  
                &W_1 \ar[u]\ar[r] 
                &W_0=\KKr \ar[u] \\
}
\]
with $\holim W_n \he *$.
The spectral sequence for this tower is isomorphic to the spectral
sequence for our Postnikov tower: in fact, there is a map of towers
$\Omega P_*(\KKr) \ra W_*$ which induces weak equivalences on the
fibers.

At this point the goal becomes to produce pairings $W_m \Smash W_n \ra
W_{n+m}$ which commute on-the-nose with the maps in the towers, and
where $W_0 \Smash W_0 \ra W_0$ is the usual multiplication on
$\Z\times BU$.  It is easy to produce pairings which commute {\it up
to homotopy\/} with the maps in the towers, and then an obstruction
theory argument shows that the maps can be rigidified.  This is the
standard argument, and has been written up in detail in \cite{D2}.
The only thing which requires much thought in the present context is
carrying out the relevant {\it equivariant\/} obstruction theory, but
this is not hard in the end.  We omit the details because the
argument is not particularly revealing.

The pairings $W_m \Smash W_n \ra W_{m+n}$ now induce a multiplicative
structure on the homotopy spectral sequence in the usual way; the
reader is again referred to \cite{D2}.

\subsection{The weight filtration and the $\gamma$-filtration}\mbox{}\par

The {\it weight filtration\/} on $KR^0(X)=[X,\KKr]$ is the one defined
by the above tower: $F^n KR^0(X)$ is defined to be the image of
$[X,W_n]$ in $[X,\KKr]$, or equivalently as the subgroup of $[X,\KKr]$
consisting of all elements which map to $0$ in $[X,P_{2n-2}(\KKr)]$.
This is a multiplicative filtration, and by Proposition~\ref{pr:adams}
it has the property that if $x\in F^nKR^0(X)$ then $\psi^k x=k^n x
\,(\text{mod}\ F^{n+1})$.  If $X$ is a space for which $H^{2n,n}(X)=0$
for $n \ggg 0$ we know that $\psi^k$ acts diagonally on
$KR^0(X)\tens\Q$, with eigenvalues $k^0$, $k^1$, $k^2$, etc.  The
tower shows that $F^n\tens \Q$ coincides with the sum of the
eigenspaces corresponding to $k^i$, for $i\geq n$.

In SGA6, Grothendieck introduced the {\it $\gamma$-filtration} on
algebraic $K^0$, designed to be an algebraic substitute for the
topological filtration induced by the classical Atiyah-Hirzebruch
spectral sequence.  For any rank $0$ stable Real bundle $\xi$ on a
$\Zt$-space $X$, one has elements $\gamma^i(\xi) \in KR^0(X)$ (see
\cite[Sec. 14]{Gr} for a nice exposition).  The $\gamma$-filtration is
defined by letting $F^n_\gamma$ be the subgroup of $KR^0(X)$ generated
by all products $\gamma^{i_1}(\xi_1)\gamma^{i_2}(\xi_2)\cdots
\gamma^{i_k}(\xi_k)$, where the $\xi_i$'s are rank 0 stable bundles
over $X$ and $i_1+i_2+\ldots+i_k \geq n$.  (So it is the smallest
multiplicative filtration in which $\gamma^i(\xi)$ is in $F^i$).  By
playing around with the algebraic definitions of $\gamma^i$ and
$\psi^k$, one can see that $F^n_\gamma \tens \Q$ also coincides with
the sum of eigenspaces of $\psi^k$ for the eigenvalues $k^i$, $i\geq
n$ (an explanation can be found in \cite[Sec. 14]{Gr}).

\begin{prop} For any $\Zt$-space $X$ one has
$F^n_\gamma KR^0(X) \subseteq F^nKR^0(X)$.  If $H^{2n,n}(X)=0$ for
$n\ggg 0$, this becomes an
equality after tensoring with $\Q$.
\end{prop}

\begin{proof}
We have already discussed the agreement rationally, since both
filtrations give eigenspace decompositions for the Adams operations.
To understand the integral story, one regards $\gamma^i$ as a map $BU
\ra \Z\times BU$ (or as an element of $KR^0(BU)$).  If $E$ is a Real
bundle of dimension $i<n$ then one can see algebraically that
$\gamma^n(E-i)=0$.  So $\gamma^n$ is null on $BU(n-1)$, and hence
factors through the homotopy cofiber $BU/BU(n-1)$.  Both $BU$ and
$BU(n-1)$ are weakly equivalent to Grassmannians, and one finds that
$\HZ^{*,*}(BU)=\HZ^{*,*}(\pt)[c_1,c_2,\ldots ]$ (where $c_i$ has
degree $(2i,i)$), and
$\HZ^{*,*}(BU(n-1))=\HZ^{*,*}[c_1,c_2,\ldots,c_{n-1}]$.  So it follows
that $\HZ^{*,*}(BU/BU(n-1))=\HZ^{*,*}(\pt)[c_n,c_{n+1},\ldots]$.  In
particular, $\HZ^{2i,i}(BU/BU(n-1))=0$ for $i< n$.  Therefore the map
$\gamma^n \colon BU/BU(n-1) \ra \Z\times BU$ lifts to $W_n$ in the
tower, and any element $\gamma^n(E-i)$ belongs to
$F^nKR^0(X)$.
\end{proof}
\vf


\section{Connective {$KR$}-theory}

The final task is to stabilize the spectral sequence we produced in
the previous section.  That spectral sequence converged to
$KR^{p+q}(X)$ only for $p+q < 0$, and we'd like to repair this
deficiency.  This is not at all difficult, and proceeds exactly as in
the non-equivariant case.  What we will do is construct a
``connective'' version of $KR$-theory, represented by a spectrum we'll
call $kr$.  There will be a homotopy cofiber sequence
\[ \Sigma^{2,1} kr \llra{\beta} kr \lra H\Zc, \]
and the Bockstein spectral sequence associated with the map $\beta$
will give the stabilized version of the spectral sequence we've been
considering.

\medskip

As in section~\ref{se:mult}, $W_n$ denotes the homotopy
fiber of $\KKr \ra \PP_{2n-2}(\KKr)$.

\begin{prop}
\label{pr:krspec}
There are weak equivalences $W_n \ra \Omega^{2,1} W_{n+1}$, unique up
to homotopy, which commute with the Bott map in the following diagram:
\[ \xymatrix{ W_n \ar[r]\ar[d] & \Omega^{2,1}W_{n+1}\ar[d] \\
              \KKr \ar[r] & \Omega^{2,1}(\KKr). 
}\]
\end{prop}

\begin{proof}[Proof of Proposition~\ref{pr:krspec}]
Consider the natural map $\alpha:\KKr \ra \PP_{2n}(\KKr)$, and apply
$\Omega^{2,1}(\blank)$.  The fact that $X\in \cA_{(2n-2,n-1)}
\Rightarrow S^{2,1}\Smash X\in \cA_{(2n,n)}$ shows that
$\Omega^{2,1}\PP_{2n}(\KKr)$ is $\cA_{(2n-2,n-1)}$-null.  By
Proposition~\ref{pr:locohoco}(a) this implies there is a lift
\[ \xymatrix{\Omega^{2,1}(\KKr) \ar[r]\ar@{ >->}[d] 
                   & \Omega^{2,1}\PP_{2n}(\KKr) \\
              \PP_{2n-2}(\Omega^{2,1}(\KKr))\ar@{.>}[ur]_{l}
}
\]
and this lift is unique up to homotopy.

Now let $\beta:\KKr\ra \Omega^{2,1}(\KKr)$ be the Bott map, and
consider the diagram
\[ \xymatrix{
   W_n \ar[r] & \KKr \ar[r]\ar[d]^{\beta} &\PP_{2n-2}(\KKr)\ar[d]^{P\beta} \\
              &\Omega^{2,1}(\KKr) \ar[r]\ar[d]^{id} 
                     &\PP_{2n-2}(\Omega^{2,1}(\KKr))\ar[d]^{l}\\ 
   \Omega^{2,1}W_{n+1}\ar[r]
              &\Omega^{2,1}(\KKr) \ar[r]
                     &\Omega^{2,1}\PP_{2n}(\KKr)
}
\]
It follows that there is a map on the homotopy fibers $W_n \ra
\Omega^{2,1}W_{n+1}$ making the diagram commute.
We need to show that this is a weak equivalence, and the procedure  is
one which should be familiar by now: we use Lemma~\ref{le:waner}.

Using Proposition~\ref{pr:pkprop}(a,b) and
Lemma~\ref{le:tildePBU}, one shows that
\smallskip
\begin{itemize}
\item  $[S^{k,0},W_n]_*=0$ for $0\leq k <n$, 
\item  $[\Zt_+\Smash S^{k,0},W_n]_*=0$ for $0\leq k<2n$, and 
\item  the same is true with $W_n$ replaced by $\Omega^{2,1}W_{n+1}$.  
\end{itemize}
The definition of $\PP_{2n-2}$ yields that the maps
\[ [S^{2n+k,n},W_n]_*\ra [S^{2n+k,n},\KKr]_* \quad \text{and} \]
\[ 
   [\Zt_+\Smash S^{2n+k,n}, W_n]_*\ra [\Zt_+\Smash S^{2n+k,n},\KKr]_*
\]
are isomorphisms for $k\geq 0$, using the fact that 
$[S^{2n+k,n},\PP_{2n-2}(\KKr)]=0$, etc.
Then the square
\[  \xymatrix{ W_n \ar[r]\ar[d] &\KKr \ar[d]^{\beta}_{\he} \\
              \Omega^{2,1}W_{n+1} \ar[r] & \Omega^{2,1}(\KKr)
}
\]
shows at once that $W_n\ra \Omega^{2,1}W_{n+1}$ induces an
isomorphism on $[S^{2n+k,n},\blank]_*$ and $[\Zt_+\Smash
S^{2n+k,n},\blank]_*$ for $k\geq 0$.  By Lemma~\ref{le:waner}, $W_n\ra
\Omega^{2,1}W_{n+1}$ is an equivalence.  
\end{proof}

\begin{defn}
Let $kr$ be the equivariant spectrum consisting of the spaces
$\{W_n\}$ and the maps $W_n\ra \Omega^{2,1}W_{n+1}$ given by the above
proposition.  The object $kr$ is called the \dfn{connective
\mdfn{$KR$}-spectrum}.
\end{defn}

The $\Omega$-spectrum for $\Sigma^{2,1}kr$ has $n$th space equal to
$W_{n+1}$, so the maps $W_{n+1} \ra W_n$ give a
`Bott map' $\Sigma^{2,1}kr \ra kr$. 
Corollary~\ref{co:ptower} identifies the homotopy fiber as 
$\Sigma^{-1,0}H\Zc$, which is equivalent to the homotopy cofiber being
$H\Zc$.  So we may form the tower of homotopy cofiber sequences
\[ \xymatrix{
\cdots \ar[r] & \Sigma^{2,1}kr\ar[d] \ar[r]^\beta 
              &kr \ar[r]^-{\Sigma^{-2,-1}\beta}\ar[d] 
              &\Sigma^{-2,-1}kr \ar[r]\ar[d] &\cdots \\
&\Sigma^{2,1}H\Zc & H\Zc & \Sigma^{-2,-1}H\Zc
}
\]
The colimit of the spectra in the tower is clearly $KR$, and the
homotopy inverse limit is contractible (these follow from thinking
about the spaces in the $\Omega$-spectra for everything in the tower).
This gives a stable version of the spectral sequence we've been
considering: for any space $X$ we have $H^{p,-\frac{q}{2}}(X)
\Rightarrow KR^{p+q,0}(X)$.  It converges conditionally because the
holim of the tower is contractible, and if $RE_\infty=0$ then it
converges strongly by \cite[Thm 8.10]{Bd} (in the language of that
result, the condition `$W=0$' is easily checked to hold).

\begin{remark}
The Postnikov tower we've constructed---and its resulting spectral
sequence---can be used to completely determine the homotopy groups of
the spaces $\PP_n(\Z\times BU)$, and hence of $W_n$ as well.  In other
words, we can completely determine the groups $kr^{*,*}(\pt)$, and in
fact the ring structure can also be deduced.  At the moment, however,
the answer doesn't seem to admit a simple description---in this sense
it is somewhat like the ring $H^{*,*}(\pt;\Zc)$, only more
complicated.  It is {\it not} true that $kr^{*,*}(\pt)\cong
H^{*,*}(\pt)[v]$, as one might naively guess based on the
non-equivariant case.  The reason essentially comes down to the fact
that there are non-trivial differentials in the spectral sequence when
$X=\pt$ (see below).  The paper \cite{HK1} computes the much more
complicated ring $M\R^{*,*}(\pt)$, and their methods can be used to
give $kr^{*,*}(\pt)$ as well.
\end{remark}

\subsection{The spectral sequence for $X=\pt$}
\mbox{}\par
\label{se:spseq4pt}

In the following diagram we draw the spectral sequence
\[ H^p(\pt;\Z(-\textstyle{\frac{q}{2}}))\Rightarrow
KR^{p+q}(\pt)=KO^{p+q}(\pt), 
\]
but using Adams indexing rather than the usual Serre conventions.
In spot $(a,b)$ we have drawn $H^{b,\frac{a+b}{2}}(\pt)$, and the
vertical line $a=N$ gives the associated graded of
$KO^{-N}$.  Said differently, the $a$-axis measures
$-(p+q)$ and the $b$-axis measures $p$.

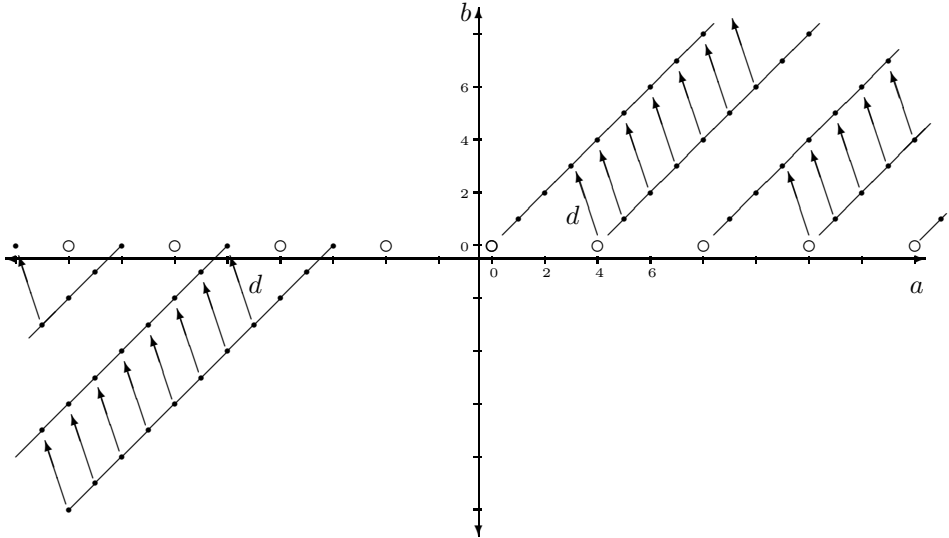
\begin{picture}(348,210)(-164,-100)
\put(-174,5){\vector(1,0){348}}
\put(174,5){\vector(-1,0){348}}
\put(5,-100){\vector(0,1){200}}
\put(5,100){\vector(0,-1){200}}
\multiput(10,10)(40,0){5}{\circle{4}}
\multiput(20,20)(10,10){8}{\circle*{2}}
\multiput(60,20)(10,10){8}{\circle*{2}}
\multiput(100,20)(10,10){7}{\circle*{2}}
\multiput(140,20)(10,10){4}{\circle*{2}}
\multiput(180,20)(10,10){1}{\circle*{2}}
\multiput(14,14)(40,0){2}{\line(1,1){80}}
\put(94,14){\line(1,1){70}}
\put(134,14){\line(1,1){42}}
\put(172,12){\line(1,1){10}}
\multiput(10,10)(-40,0){5}{\circle{4}}
\multiput(-50,10)(-10,-10){11}{\circle*{2}}
\multiput(-90,10)(-10,-10){8}{\circle*{2}}
\multiput(-130,10)(-10,-10){4}{\circle*{2}}
\multiput(-170,10)(-10,-10){1}{\circle*{2}}
\put(-50,10){\line(-1,-1){100}}
\put(-90,10){\line(-1,-1){80}}
\put(-130,10){\line(-1,-1){35}}
\multiput(59,22)(10,10){6}{\vector(-1,3){8}}
\multiput(139,22)(10,10){4}{\vector(-1,3){8}}
\multiput(50,14)(80,0){2}{\vector(-1,3){8}}
\multiput(-81,-18)(-10,-10){8}{\vector(-1,3){8}}
\multiput(-161,-18)(-10,-10){1}{\vector(-1,3){8}}
\multiput(-170,3)(20,0){18}{\line(0,1){3}}
\multiput(3,-90)(0,20){10}{\line(1,0){3}}
\put(38,18){$d$}
\put(-82,-8){$d$}
\put(9,-2){$\scriptscriptstyle{0}$}
\put(29,-2){$\scriptscriptstyle{2}$}
\put(49,-2){$\scriptscriptstyle{4}$}
\put(69,-2){$\scriptscriptstyle{6}$}
%
\put(-2,8){$\scriptscriptstyle{0}$}
\put(-2,28){$\scriptscriptstyle{2}$}
\put(-2,48){$\scriptscriptstyle{4}$}
\put(-2,68){$\scriptscriptstyle{6}$}
\put(-2,95){${b}$}
\put(168,-8){${a}$}
\end{picture}

There are several points to make:

\begin{enumerate}[(a)]
\item Using the multiplicative properties of the spectral sequence,
one only has to determine the two differentials labelled `$d$'---all
the others can be deduced from these.
Since we know the groups $KO^*(\pt)$, it's clear that these two
differentials have to exist. (It would be nice to have a more
intrinsic explanation, however).  
\item The spectral sequence collapses at the next page.
\item The unstable spectral sequence of sections 4 and 5 is the part
in the first quadrant.  We can read off the action of the Adams
operations on $KO^n(\pt)$ for $n\leq 0$ directly from the diagram, the
`weight lines' being along the antidiagonals: $KO^0$ is pure of weight
$0$, $KO^{-1}$ is pure of weight $1$, $KO^{-2}$ and $KO^{-4}$ are of
weight $2$, $KO^{-8}$ is of weight $4$, etc.
\item 
Everything about the ring structure on $KO^*$, as well as the
comparison map $KO^* \ra K^*$, can be read off of this spectral
sequence and the corresponding spectral sequence where $X=\Zt$ (which
lies entirely along the $b=0$ line, and hence collapses).  They can be
deduced from our knowledge of $H^{*,*}(\pt)$ and $H^{*,*}(\Zt)$
provided by Theorem~\ref{th:coeffs}.
\item The part of the spectral sequence in the first quadrant is known
to topologists in another setting: it's the Adams spectral sequence
for $bo$ based on $bu$. 
\end{enumerate}

\vf

\section{\'Etale analogs}
\label{se:etale}
The difference between algebraic $K$-theory and \'etale $K$-theory, or
motivic cohomology and \'etale motivic cohomology, is very familiar in
the motivic setting.  In this section we play with similar ideas in
the $\Zt$ world.  The analogs are well known, although the only source
seems to be \cite[Section 3.3]{MV}, which doesn't develop things in
much detail.  We use these ideas to give a proof of the classical fact
that $K^{\Zt} \he KO$ and $(\Z\times BU)^{h\Zt} \he \Z\times BO$.

\medskip

Let us return briefly to the setting where $G$ is any finite group.
By an \dfn{equivariant covering space} $E\ra B$ we mean an equivariant
map which, after forgetting the $G$-actions, is a covering space in
the usual sense.  Given such a map we may form its \CCech complex 
$\Cech(E)$, which is the simplicial space
\[ \xymatrix{ E &
   E\times_B E \ar@<0.5ex>[l]\ar@<-0.5ex>[l] & E\times_B E\times_B E
     \ar@<0.6ex>[l]\ar[l]\ar@<-0.6ex>[l] 
    & \cdots
       \ar@<0.6ex>[l]\ar@<0.2ex>[l]\ar@<-0.2ex>[l]\ar@<-0.6ex>[l]   
}
\] 
(where we have omitted the degeneracies for typographical reasons).
In non-equivariant topology, the map $\hocolim_n \Cech(E)_n \ra B$ is
a weak equivalence (cf. \cite[Cor. 1.3]{DI}).  This is not true
equivariantly, as the covering space $G\ra *$ shows.  In this case the
realization of the \CCech complex is precisely $EG$, and $EG\ra *$ is
not an equivariant  equivalence.  We will see that in some sense this
turns out to be the only problem, though.

If $Z$ is a $G$-space there is a natural map of $G$-spaces
\[ 
F(B,Z) \ra \holim_n F(\Cech(E)_n,Z)
\]
(here $F(X,Y)$ is the usual mapping space, with its induced $G$-action).
We will say the space $Z$ satisfies \mdfn{\'etale descent for the
covering $E\ra B$} if this natural map is an equivariant weak
equivalence.  This is the same as requiring that the corresponding
maps for the coverings $G/H \times E \ra G/H \times B$ all be
nonequivariant equivalences, where $H$ ranges over the subgroups of
$G$.  If the phrase is not qualified, then ``\'etale descent'' means
``\'etale descent for all covering spaces''.

Using results of \cite[Chaps. 3,4]{H}, we may localize the model
category $\Top_G$ at the maps $\hocolim \Cech(E) \ra B$ where $E\ra B$
ranges over the elements of the set $\{G/H \times G \ra G/H \times {*}
| H < G\}$.  This produces a new model category structure which we'll
denote $\Top_G^{et}$.  As the following result shows, the fibrant
objects are precisely the spaces which satisfy \'etale descent for all
covering spaces.


\begin{prop}\mbox{}\par
\begin{enumerate}[(a)]
\item If $Z$ is any $G$-space, then $Z^{EG}$ satisfies \'etale descent.
\item A map $X\ra Y$ is a weak equivalence in $\Top_G^{et}$ iff it is
a non-equivariant weak equivalence.
\item For any space $Z$, the map $Z \ra Z^{EG}$ is a fibrant
replacement in $\Top_G^{et}$.  
\end{enumerate}
\end{prop}

\begin{proof}
Let $X\ra Y$ be a $G$-map which is also a nonequivariant equivalence,
and assume that $X$ and $Y$ are cofibrant.  Then $X\times EG \ra
Y\times EG$ is an equivariant equivalence, and therefore so is
$F(Y\times EG,Z) \ra F(X\times EG,Z)$ for any $Z$.  By adjointness
this map is the same as $F(Y,Z^{EG}) \ra F(X,Z^{EG})$.  In particular,
if $E\ra B$ is an equivariant covering space then by taking $X\ra Y$
to be the map $\hocolim \Cech(E) \ra B$ (which is a nonequivariant
equivalence by \cite[Cor. 1.3]{DI}) we find that $Z^{EG}$ satisfies
\'etale descent.  This proves (a).

When forming $\Top_G^{et}$ we are localizing at maps which are
nonequivariant equivalences.  It follows from this that every
equivalence in $\Top_G^{et}$ is a nonequivariant equivalence.  For
(b), we must show the other direction.  Note that if $Z$ is a fibrant
object in $\Top_G^{et}$, then $Z \ra Z^{EG}$ is an equivariant
equivalence (this is \'etale descent for $G\ra *$).  If $X\ra Y$ is a
map between cofibrant objects, we may consider the diagram
\[ \xymatrix{ 
F(Y,Z) \ar[d]\ar[r]^-\sim & F(Y,Z^{EG}) \ar[d]\ar[r]^-\cong 
                      & F(Y\times EG,Z)\ar[d] \\
F(X,Z) \ar[r]^-\sim & F(X,Z^{EG}) \ar[r]^-\cong & F(X\times EG,Z)}.
\]
If $X\ra Y$ was a nonequivariant equivalence then $X\times EG \ra
Y\times EG$ is an equivariant equivalence, and so the right vertical
map is an equivalence as well.  It follows that $F(Y,Z)\ra F(X,Z)$ is
an equivariant equivalence for every fibrant object $Z$ in
$\Top_G^{et}$, and therefore $X\ra Y$ is an equivalence in
$\Top_G^{et}$.

Part (c) is an immediate consequence of (a) and (b).
\end{proof}

We will call $Z^{EG}$ the \dfn{\'etale localization} (or {\it Borel\/}
localization) of the space $Z$, and we'll sometimes write it as
$Z_{et}$.  Note that \'etale localization preserves fiber sequences and
homotopy limits.

Everything from our above discussion generalizes directly to spectra
as well.  We can talk about an equivariant spectrum $\cE$ which
satisfies \'etale descent---these correspond to what are usually
called `Borel cohomology theories' \cite[p. 233]{M}.  If $\cE$ is the
$RO(G)$-graded spectrum given by $V\ra \cE_V$, its \'etale
localization (or corresponding Borel theory) is the spectrum
$\cE_{et}$ given by $V\ra \cE_V^{EG}$.  Note that if $\cE$ was an
$\Omega$-spectrum then $\cE_{et}$ is also an $\Omega$-spectrum.

\medskip

At this point we switch back to the $\Zt$-setting, where we
can write down the following two results.  The first
is an immediate consequence of Corollary~\ref{co:ptower}, the second
of Theorem~\ref{th:stabmain}.

\begin{prop}  
\label{co:ptoweret}
There is a tower of homotopy fiber sequences
\[
   \xymatrix{   
                \quad \cdots\quad \ar[r]
                &[\PP_{4}(\KK)]_{et} \ar[r]  
                &[\PP_{2}(\KK)]_{et} \ar[r] 
                &[\PP_{0}(\KK)]_{et} \ar[r]
                &{*}\\            
            &K(\Z(2),4)_{et}\ar[u]
            &K(\Z(1),2)_{et}\ar[u]
            &\Z_{et} \ar[u]
}
\]
and the homotopy limit of the tower is $(\KK)_{et}$.
\end{prop}

\begin{prop}
\label{pr:ettowerstab}
There is a tower of homotopy cofiber sequences in $\Zt$-spectra of the form
\[ \xymatrix{
\cdots \ar[r] & \Sigma^{2,1}kr_{et}\ar[d] \ar[r]^\beta 
              &kr_{et} \ar[r]^-{\Sigma^{-2,-1}\beta}\ar[d] 
              &\Sigma^{-2,-1}kr_{et} \ar[r]\ar[d] &\cdots \\
&\Sigma^{2,1}H\Zc_{et} & H\Zc_{et} & \Sigma^{-2,-1}H\Zc_{et}
}
\]
The homotopy colimit of the tower is $KR_{et}$, and the homotopy limit
is contractible.
\end{prop}

We want to analyze the spectral sequence which comes from the above
tower, so to start with we need a knowledge of $H\Zc_{et}$.  
The theory $H\Zc_{et}$ has coefficient groups described as follows:

\begin{prop}\mbox{}\par
\label{pr:HZandHZet}
\begin{enumerate}[(a)]
\item $\Het^{*,*}(\pt)=\Z[x,x^{-1},y]/(2y)$ where $x$ has degree
$(0,2)$ and $y$ has degree $(1,1)$.  
\item $\Het^{*,*}(\Zt)=\Z[u,u^{-1}]$  where $u$ has degree $(0,1)$.
\item The Mackey functor $\underline{H\Z}_{et}^{0,2n}$ is equal to
$\Zm$, for all $n$.
\item The map $\HZ^{*,*}(\Zt) \ra \Het^{*,*}(\Zt)$ is an isomorphism.
The map $\HZ^{p,q}(\pt) \ra \Het^{p,q}(\pt)$ is an isomorphism when
$p<2q$, and multiplication by $2$ when $p=0$ and $q<0$ is even.  These are
the only degrees in which its possible to have a nonzero map. 
\end{enumerate}
\end{prop}

The proof is given in Appendix B.

Proposition~\ref{pr:ettowerstab} gives rise to a spectral sequence
$\Het^{p,-\frac{q}{2}}(X) \Rightarrow KR_{et}^{p+q,0}(X)$.   The
spectral sequence for $X=\pt$ is drawn below, using the same indexing
conventions as in section~\ref{se:spseq4pt}.

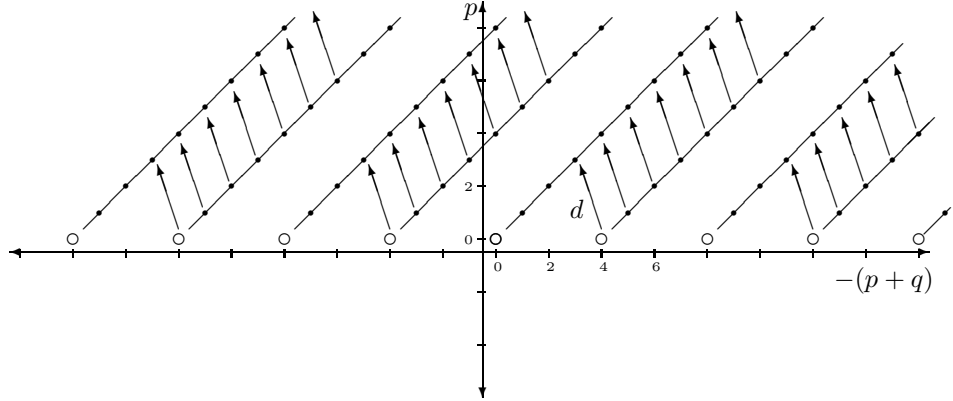
\begin{picture}(348,160)(-164,-50)
\put(-174,5){\vector(1,0){348}}
\put(174,5){\vector(-1,0){348}}
\put(5,-50){\vector(0,1){150}}
\put(5,100){\vector(0,-1){150}}
%
\multiput(10,10)(40,0){5}{\circle{4}}
\multiput(20,20)(10,10){8}{\circle*{2}}
\multiput(60,20)(10,10){8}{\circle*{2}}
\multiput(100,20)(10,10){7}{\circle*{2}}
\multiput(140,20)(10,10){4}{\circle*{2}}
\multiput(180,20)(10,10){1}{\circle*{2}}
\multiput(14,14)(40,0){2}{\line(1,1){80}}
\put(94,14){\line(1,1){70}}
\put(134,14){\line(1,1){42}}
\put(172,12){\line(1,1){10}}
\multiput(59,22)(10,10){6}{\vector(-1,3){8}}
\multiput(139,22)(10,10){4}{\vector(-1,3){8}}
\multiput(50,14)(80,0){2}{\vector(-1,3){8}}
%
\multiput(10,10)(-40,0){5}{\circle{4}}
\multiput(-20,20)(10,10){8}{\circle*{2}}
\multiput(-60,20)(10,10){8}{\circle*{2}}
\multiput(-100,20)(10,10){8}{\circle*{2}}
\multiput(-140,20)(10,10){8}{\circle*{2}}
\multiput(-26,14)(-40,0){4}{\line(1,1){80}}
%
\multiput(-21,22)(10,10){6}{\vector(-1,3){8}}
\multiput(-101,22)(10,10){6}{\vector(-1,3){8}}
\multiput(-30,14)(-80,0){2}{\vector(-1,3){8}}

\put(38,18){$d$}
%
\multiput(-170,3)(20,0){18}{\line(0,1){3}}
\multiput(3,-30)(0,20){7}{\line(1,0){3}}
\put(9,-2){$\scriptscriptstyle{0}$}
\put(29,-2){$\scriptscriptstyle{2}$}
\put(49,-2){$\scriptscriptstyle{4}$}
\put(69,-2){$\scriptscriptstyle{6}$}
\put(-2,8){$\scriptscriptstyle{0}$}
\put(-2,28){$\scriptscriptstyle{2}$}
\put(-2,95){${p}$}
\put(138,-8){${-(p+q)}$}
\end{picture}

Note that there is a map from the spectral sequence of
Theorem~\ref{th:stabspseq} to the one above (coming from the maps $\HZ
\ra \Het$ and $KR \ra KR_{et}$).  Based on our discussion in
section~\ref{se:spseq4pt} we therefore know that the differential
labelled `$d$' has to exist---in fact, we know all the differentials
in the first quadrant below the $y=x$ line.  Multiplicativity then
allows us to deduce what's happening in the rest of the spectral
sequence, and that is what is shown above.  As we see, the spectral
sequence again converges to $KO^*$; but this time we've actually
gained some information:

\begin{cor}
\label{co:KRet}
The natural map $KR \ra KR_{et}$ is a weak equivalence---that is, $KR$
satisfies \'etale descent.
\end{cor}

\noindent
In very fancy language, this says `the Quillen-Lichtenbaum
conjecture holds for $KR$'.

\begin{proof}
The map is of course a nonequivariant equivalence, so we only have to
analyze what happens on fixed sets---i.e., we study the map
$\alpha\colon KR^{*,0}(\pt) \ra KR_{et}^{*,0}(\pt)$. Consider the map
of spectral sequences from $H\Z^{p,-\frac{q}{2}}(\pt) \Rightarrow
KR^{p+q,0}(\pt)$ to the corresponding \'etale version.  Using
everything we know about both spectral sequences, as well as what
Proposition~\ref{pr:HZandHZet}(d) says about the map
$\HZ^{*,*}(\pt) \ra \Het^{*,*}(\pt)$, it is easy to see that $\alpha$
is an isomorphism when $*\leq 0$, and also when $*=-4n$.  The fact
that $\alpha$ is a ring map then gives us the isomorphism in the
remaining dimensions.
\end{proof}

\begin{cor}
\label{co:BUdesc}
If we consider $\Z\times BU$ as a nonequivariant space with its $\Zt$
action, then $(\Z\times BU)^{h\Zt}\he \Z\times BO$.  Similarly,
$K^{h\Zt}\he KO$.  
\end{cor}

\begin{proof}
This is just a translation of the previous corollary.  The $0$th space
in the $\Omega$-spectrum of $KR$ is $\Z\times BU$, whereas the $0$th
space in the $\Omega$-spectrum of $KR_{et}$ is $(\Z\times BU)^{E\Zt}$.
Corollary~\ref{co:KRet} tells us that these spaces are equivariantly
equivalent, and therefore they have weakly equivalent fixed sets.  
The proof that $K^{h\Zt}\he KO$ follows the
same lines, but takes place in the stable category.
\end{proof}

Corollary~\ref{co:BUdesc} (and \ref{co:KRet}, which is equivalent) is
of course well known---a recent reference is \cite{K}.  In the end our
proof is only slightly different from the classical proof which
analyzes $\pi_*K^{h\Zt}$ via the spectral sequence for $\pi_*$ of a
homotopy limit, making use of the map $KO\ra K^{h\Zt}$ (in fact this
spectral sequence has the same form as the one drawn above).  It may
be worth summarizing our proof to exhibit the similarities.  The
spectral sequence relating $H\Z$ and $KR$, when applied to a point,
gave us something converging to $KO$ (the fixed set of $KR$).  Because
we knew the homotopy groups of $KO$, we could analyze the
differentials in this spectral sequence.  On the other hand we have a
corresponding spectral sequence relating $H\Z_{et}$ and $KR_{et}$;
when applied to a point, it converges to $K^{h\Zt}$ (the fixed set of
$KR_{et}$).  There is a comparison map of spectral sequences, and our
knowledge of the differentials in the first lets us deduce the
differentials in the second.

\vf

\section{Postnikov sections of spheres}
\label{se:symmprod}

In this final section we prove Theorem~\ref{th:P(S)}, which identifies
the Postnikov section $\PP_V(S^V)$ with the Eilenberg-MacLane space
$K(\Zm,V)$ over the group $G=\Z/2$.  Theorem~\ref{th:main}(a), which
is all we really need in this paper, is the case where $V=\C^n$.  This
section takes place entirely in the context of $\Z/2$-spaces.

\bigskip

Suppose that $V$ contains a copy of the trivial representation.  It
follows from \cite{dS} and Corollary~\ref{co:symm=ag} of the present
paper that the infinite symmetric product $\SP^\infty(S^V)$ is a model
for $K(\Zm,V)$ (this actually works over any finite group, not just
$G=\Z/2$).  Using this, we have:

\begin{lemma}
\label{le:locsphere}
If $V\supseteq 1$ then $\Spi(S^{V}) \ra \PP_{V}(\SPi(S^{V}))$ is a
weak equivalence.
\end{lemma}

\begin{proof}
We know by Theorem~\ref{th:coeffs} that for any $r,s \geq 0$
\begin{eqnarray*}
{[S^{V+r\R+s\R_{-}},K(\Zm,V)]_*} &=& H^{-r-s,-s}(\pt;\Zc) = 0 \quad 
\text{and}, \\
 {[\Zt_+\Smash S^{V+r\R+s\R_{-}},K(\Zm,V)]_*} &=& H^{-r-s,-s}(\Zt;\Zc) = 0 
\end{eqnarray*}
as long as $r$ and $s$ are not both zero.  This shows that
$\Spi(S^{V})$ is $\cA_{V}$-null, which implies that the map in the
statement of the lemma is a weak equivalence.
\end{proof}

Our goal is the following:

\begin{thm}
\label{th:locsphere}
If $V \supseteq 1$ and $S^V \ra \SPi(S^V)$ is the obvious map, then
\[ \PP_{V}(S^{V}) \ra \PP_{V}(\Spi(S^{V}))
\]
is a weak equivalence.  Therefore $\PP_{V}(S^{V})\he
K(\Z,V)$.
\end{thm}

\noindent
Note that the second statement follows immediately from the first,
in light of the above lemma.  The proof of the first statement is
based on a geometric analysis of the infinite symmetric product
construction, and involves the following steps:

\smallskip

\begin{enumerate}[(1)]
\item For each $k\geq 2$ we produce a homotopy cofiber sequence
\[ S^{V}\Smash \Bigl([V\tens \tilde{R}]-0\Bigr )/\Sigma_k 
\ra \SP^{k-1}(S^{V}) \ra
    \SP^k(S^{V})
\]
where $\tilde{R}$ is the reduced standard representation of the
symmetric group $\Sigma_k$ (defined below) equipped with the trivial
$\Z/2$ action.
\item We show that 
\[ \PP_{V} \Bigl (S^{V}\Smash ([V\tens \tilde{R}]-0)/\Sigma_k \Bigr) 
\he *.\]
The key ingredient for this is a geometric analysis of the 
$\Zt$--fixed sets of the orbit space $([V\tens \tilde{R}]-0)/\Sigma_k$.
\item From (1) and (2) it follows that 
\[ \PP_{V}(\SP^{k-1}(S^{V})) \ra \PP_{V}(\SP^k(S^{V}))\]
is a weak equivalence for every $k\geq 2$, and passing to the limit
yields the same for $\PP_{V}(S^{V}) \ra
\PP_{V}(\SPi(S^{V}))$.
\end{enumerate}

\medskip

\subsection{Step 1}

Consider the filtration of $\SPi(S^{V})$ given by the finite
symmetric products:
\[ S^{V} \subseteq \SP^2(S^{V}) \subseteq SP^3(S^{V})
\subseteq \cdots \subseteq \Spi(S^{V}).\]
Recall that the \dfn{standard representation} of the symmetric group
$\Sigma_k$ is the space $R=\R^k$ where the group acts by permuting the
standard basis elements.  This contains a trivial, one-dimensional
subrepresentation consisting of all vectors $(r,r,\ldots,r)$, and the
\dfn{reduced standard representation} $\tilde R$ is the quotient of $R$
by this subrepresentation.  We regard $R$ and $\tilde{R}$ as
$\Z/2$-representations by giving them the trivial actions.

The following proposition was inspired by \cite[Thm. 2.3]{JTTW},
which handled the case $k=2$:

\begin{prop}
\label{pr:spdecomp}
The inclusion $\SP^{k-1}(S^{V}) \inc \SP^k(S^{V})$ sits in a
homotopy cofiber sequence of the form
\[ S^{V} \Smash ([V\tens \tilde{R}]-0)/\Sigma_k \ra SP^{k-1}(S^{V}) \ra
\SP^k (S^{V}) \]
where $\tilde R$ denotes the reduced standard representation of $\Sigma_k$.
\end{prop}

\begin{proof}
To save ink, write $B=B(V)$ and $S=S(V)$ for the unit
ball and unit sphere in $V$.  
We begin with the relative homeomorphism
\[ (B,S) \mapiso (S^{V},*).\]
Applying $\SP^k$ to these pairs gives a relative homeomorphism
\[ (\SP^k(B), Z/\Sigma_k) \mapiso (\SP^k(S^{V}), \SP^{k-1}(S^{V}))\]
where $Z$ is the space
\[ (S\times B\times \cdots\times B) \cup
    (B \times S \times B \times\cdots\times B)
\cup \cdots\cup
   (B\times\cdots\times B\times S)\subseteq B^k.\]
This says that there is a pushout square of the form
\[ \xymatrix{
             Z/\Sigma_k \ar[r] \ar@{ >->}[d] 
                     &\SP^{k-1}(S^{V})\ar@{ >->}[d] \\
             \SP^k(B) \ar[r] & \SP^k(S^{V}).
}\]
Since $\SP^k(B)$ is clearly contractible, the desired cofiber
sequence will follow if we can identify $Z/\Sigma_k$ with
$S^{V}\Smash ([V\tens \tilde{R}]-0)/\Sigma_k$.

$Z$ naturally includes into $(V\oplus\cdots\oplus V)-0=[V\tens R]-0$, and 
the assignment
\[
  Z \ra S(V\oplus\cdots\oplus V), \qquad
  v \ra \textstyle{\frac{v}{\mid{v}\mid}}
\]
%
is a homeomorphism which is both $\Zt$- and $\Sigma_k$- equivariant.
So we may identify the $\Zt$-spaces $Z/\Sigma_k$ and $S(V\tens
R)/\Sigma_k$.

Since $R$ decomposes as $\R\oplus \tilde{R}$ (as both $\Sigma_k$- and
$\Zt$-representations), we have the corresponding decomposition
$V\tens R \iso V \oplus [V\tens \tilde{R}]$.  
Lemma~\ref{le:spherejoin} below gives a homeomorphism
\[ S\Bigl (V\oplus [V\tens \tilde{R}] \Bigr)/\Sigma_k \cong S(V) * [S(V\tens
R)/\Sigma_k], 
\] 
where $X*Y$ denotes the usual \dfn{join} of $X$ and $Y$.
It is a general fact (true in any model category) that for pointed
spaces $X$ and $Y$, the join $X*Y$ is weakly equivalent to
$\Sigma(X\Smash Y)$.  Because
$V\supseteq 1$, both $S(V)$ and $S(V\tens \tilde{R})$ have
nonempty $\Zt$-fixed sets, and therefore can be made pointed.  So we
finally conclude that
\begin{eqnarray*}
Z/\Sigma_k \cong S(V)* [S(V\tens \tilde{R})/\Sigma_k] 
        &\he& S^1\Smash S(V) \Smash [S(V\tens \tilde{R})/\Sigma_k] \\
        &\he& S^V\Smash ([V\tens \tilde{R}]-0)/\Sigma_k.
\end{eqnarray*}
\end{proof}

\begin{lemma}
\label{le:spherejoin}
Let $V$ and $W$ be orthogonal representations of $\Zt$.  Let $G$ be a
finite group acting $\Zt$-equivariantly and orthogonally on $W$, and
let $G$ act on $V$ trivially. Then there is a natural
$\Zt$-equivariant homeomorphism
\[ S(V)*[S(W)/G] \mapiso S(V\oplus W)/G  \]
where the space on the left denotes the join.
\end{lemma}

\begin{proof}
A point in the left-hand space can be represented by a triple
$(v,t,[w])$ where $v\in S(V)$, $t\in [0,1]$, $w\in S(W)$, and $[w]$
denotes the $G$-orbit of $w$.
We leave it to the reader to check that the map
$(v,t,[w]) \assign [\sqrt{1-t}\cdot v \oplus \sqrt{t}\cdot w]$
is well-defined and a $\Zt$-equivariant homeomorphism.
\end{proof}

\subsection{Step 2}

\begin{prop}
\label{pr:lensspace}
Suppose $V=\R^p \oplus (\R_{-})^q$, where $p\geq 1$.
Let $X=([V\tens \tilde{R}]-0)/\Sigma_k$.
\begin{enumerate}[(a)]
\item The fixed set $X^{\Zt}$ is path-connected for $k\geq 3$.
\item When $k=2$, $X^{\Zt} \he \RP^{p-1} \copr \RP^{q-1}$ (where the
second summand is interpreted as $\emptyset$ when $q=0$, and a point
when $q=1$).
When 
$q\geq 1$ there
exists a map $S^{1,1} \ra X$ which on fixed sets
induces an isomorphism on $\pi_0$.
\end{enumerate}
\end{prop}

\begin{proof}
The argument involves producing explicit paths in the fixed sets.  As
it's somewhat lengthy, we postpone it until the very end of the section.
\end{proof}

\begin{cor}
\label{co:lensnull}
Consider the set $\cA$ consisting of the objects
\[ \bullet\ S^{n,0}, \ n\geq 1; \qquad\qquad
   \bullet\ S^{n,0}\Smash \Zt_+, \ n\geq 1; \qquad\qquad
   \bullet\ S^{1,1}.
\]
Then the nullification $P_{\cA} X$ at this set is contractible.
\end{cor}

\begin{proof}
Let $PX$ denote the nullification of $X$.  To show $PX$ is contractible
we have only to check that $[S^{n,0},X]_*=0=[S^{n,0}\Smash \Zt_+,X]_*$ for
all $n\geq 0$.  For $n\geq 1$ this follows just from the definition of
$PX$.  We are therefore reduced to checking $n=0$, which is the
statement that both $PX$ and $(PX)^\Zt$ are path-connected (as
non-equivariant spaces).

Clearly $X$ is path connected.  Attaching cones on maps cannot
disconnect the space, so $PX$ must also be path-connected.

Proposition~\ref{pr:lensspace} says that $X^\Zt$ is path-connected for
$k\geq 3$ (or $k=2$ and $q=0$), and attaching cones on maps cannot
disconnect the fixed set.  So $(PX)^\Zt$ is again path-connected in
this case.

When $k=2$ and $q\geq 1$ the proposition says that $X^\Zt$ has {\it
two} path components, but they are linked by an $S^{1,1}$.  Attaching
a cone on this map will give a space whose fixed set is connected, and
then reasoning as in the previous paragraph we find that $PX$ will
also have that property.
\end{proof}

\begin{cor}
The space $\PP_{V}(S^{V}\Smash ([V\tens \tilde{R}]-0)/\Sigma_k)$ is
contractible.
\end{cor}

\begin{proof}
Let $X=([V\tens \tilde{R}]-0)/\Sigma_k$, and suppose we cone off a map
$S^{1,0}\ra X$ to make a space $X'$:
\[ S^{1,0} \ra X \ra X'.\]
Smashing with $S^{V}$ gives a cofiber sequence
\[ S^{V+1} \ra S^{V}\Smash X \ra S^{V}\Smash X', \]
and since $\PP_{V}(S^{V+1})\he *$ it follows by
Proposition~\ref{pr:locohoco}(c) that 
\[ \PP_{V}(S^{V}\Smash X) \we \PP_{V}(S^{V}\Smash X').\]
In other words, we may cone off arbitrary maps $S^{1,0}\ra X$ without
effecting the homotopy type of $\PP_{V}(S^{V}\Smash X)$.
The same reasoning shows we can cone off maps $S^{1,1}\ra X$,
$S^{n,0}\ra X$, and $\Zt_+\Smash S^{n,0}\ra X$ ($n\geq 1$) with the
same result.  So the conclusion is that
\[ \PP_{V}(S^{V}\Smash X) \he \PP_{V}(S^{V}\Smash PX), \]
where $PX$ denotes the nullification considered in
Corollary~\ref{co:lensnull}. But that corollary says that $PX$ is
contractible, and so we're done.
\end{proof}

\subsection{Step 3}

\begin{proof}[Proof of Theorem~\ref{th:locsphere}]
We are to show that the map $S^{V}\ra \Spi(S^{V})$
becomes a weak equivalence after applying $\PP_{V}$.
We'll simplify
$\PP_{V}(X)$ to just $P(X)$, and $\SP^k(S^{V})$ to just $\SP^k$.
Proposition~\ref{pr:spdecomp} gives cofiber sequences
\[ S^{V}\Smash ([V\tens \tilde{R}]-0)/\Sigma_k \ra
\SP^{k-1} \ra \SP^k,
\]
and we have seen that
\[ P(S^{V}\Smash ([V\tens \tilde{R}]-0)/\Sigma_k) \he *. \]
So Proposition~\ref{pr:locohoco}(c) shows that $P(\SP^{k-1})\ra
P(\SP^k)$ is a weak equivalence.  Hence, one has a sequence of weak
equivalences
\[ P(S^{V}) \we P(\SP^2)
                       \we P(\SP^3) \we \cdots
\]
and therefore $P(S^V) \ra \hocolim_k P(SP^k)$ is a weak equivalence as
well.  

Now look at the composite of the two maps
\[ P(S^{V}) \lra P(\hocolim \SP^k) \lra
P(\colim \SP^k)=P(\Spi).\] The middle object may be identified with
$\hocolim_k P(SP^k)$ using Proposition~\ref{pr:locohoco}(b) and some
common sense, and so the first map is an equivalence.  The second map
is a weak equivalence because $\hocolim \SP^k \ra \colim \SP^k$ was
one.  Hence, the composite is also a weak equivalence.
\end{proof}

\subsection{Loose ends: The analysis of the fixed sets 
\mdfn{$ [(V\tens \tilde{R}-0)/\Sigma_k]^\Zt$}}\mbox{}\par

\smallskip

The one thing still hanging over our heads is the

\begin{proof}[Proof of Proposition~\ref{pr:lensspace}]
Recall that $X=[V\tens\tilde{R}]-0/\Sigma_k$.  Begin by decomposing $V$
as a sum of irreducibles $V=U_0 \oplus U_1 \oplus \cdots \oplus U_n$
where $U_0=1$ (or course the only irreducible representations of $\Zt$
are $\R$ and $\R_{-}$, but it's easiest to think in slightly more
generality here).  An element of $X$ is represented by a coset
$[u_0,\ldots,u_n]$ with $u_i\in U_i \tens \tilde{R}$ and at least one $u_i$
nonzero.

We begin with part (a), which says that when $k\geq 3$ the fixed set
$X^{\Zt}$ is path connected.  First note that if $u=[u_0,\ldots,u_n]\in
X^{\Zt}$ and $u_i\neq 0$, then 
\[ t \mapsto [tu_0,tu_1,\ldots,tu_{i-1},u_i,tu_{i+1},\ldots,tu_n] \]
gives a path in $X^{\Zt}$ from $u$ to $[0,0,\ldots,0,u_i,0,\ldots,0]$.
It follows that it suffices to prove the result when $V$ is of the
form $1\oplus U$, where $U$ is a (possibly $0$) irreducible
representation.  (Recall that the result requires $V$ to contain $1$,
which is why we don't reduce to $V=U$.)
Since our group is $\Zt$, we only have to worry about
$V=\R$ and $V=\R\oplus \R_{-}=\C$.  

The case $V=\R$ is trivial, because then
$X^{\Zt}=X=(\tilde{R}-0)/\Sigma_k$, and $\tilde{R}-0$ was path
connected because $\dim \tilde{R}\geq 2$ (recall $k\geq 3$ here).  So we
are left to deal with $V=\C$.  In this case $V\tens \tilde{R}$ is the
complex reduced standard representation, which we may identify with
$\{(z_1,\ldots,z_k)\in \C^k \mid \sum z_i=0\}$.  We will write
$[z_1,\ldots,z_k]$ for the coset of $(z_1,\ldots,z_k)$ in $(V\tens
\tilde{R})/\Sigma_k$.  

Let $A=\{[r_1,\ldots,r_k] \mid r_i\in \R\} \subseteq X^\Zt$.  Clearly
$A$ is path-connected, as $A\cong (\R^{k-1}-0)/\Sigma_k$.  We will
show that any element of $X^\Zt$ can be connected by a path to an
element of $A$.  If $[z_1,\ldots,z_k]\in X^\Zt$ then there is a
$\sigma\in \Sigma_k$ with the property that
\[ (z_{\sigma(1)},\ldots,z_{\sigma(k)})=(\bar{z}_1,\ldots,\bar{z}_k).
\]
By writing $\sigma$ as a composite of disjoint cyclic permutations,
it's easy to see that this can only happen if $(z_1,\ldots,z_k)$ has
the form
\[ (w_1,\bar{w}_1,\ldots,w_l,\bar{w}_l,r_1,\ldots,r_j) \]
up to permutation of the $z$'s (where $r_i\in \R$).
If all the $w_i$'s are real, then our point is already in $A$
and we can stop.  So we can assume that $w_1\not\in \R$.
Consider the path
\[ 
     t \assign
       [w_1+f(t),\bar{w}_1+f(t),tw_2,t\bar{w}_2,\ldots,tw_l,
        t\bar{w}_l,tr_1,\ldots,tr_j]
\]
where
\[ 
f(t)=-\frac{1}{2}(2\re(w_1)+2t\re(w_2)+\ldots+2t\re(w_l)+tr_1+\ldots+tr_j)
\]
(so $f(t)$ is the real number which makes the sum of the components zero
in the previous expression).
It's easy to see that this describes a path in $X^\Zt$ connecting our
original point with 
\[ 
      [w_1+f(0),\bar{w}_1+f(0),0,\ldots,0], 
\] 
and this
point has the form $[bi,-bi,0,\ldots,0]$ for some nonzero $b\in \R$.
Next we consider the path
\[ 
      t \assign [t+b(1-t)i,t-b(1-t)i,-2t,0,0,\ldots,0] 
\]
(and here we use the fact that $k\geq 3$).
This is a path in $X^\Zt$ connecting $[bi,-bi,0,\ldots,0]$ with
$[1,1,-2,0,\ldots,0]$, the latter of which is in $A$.
So this completes the proof that $X^\Zt$ is path-connected when
$V=\C$, and we are done with part (a).

\smallskip

Now we turn to part (b), which is the case $k=2$.
The reduced standard representation of $\Sigma_2$ is $\tilde R=\R$ with
the $\Sigma_2$-action equal to multiplication by $-1$.
If $V=\R^p\oplus \R_{-}^q$ then
$X=[(\R^p \oplus \R_{-}^q)-0]/\Sigma_2 \he \RP^{p+q-1}$.  Using
homogeneous coordinates $[r_1,\ldots,r_p,r_{p+1},\ldots,r_{p+q}]$ on
$\RP^{p+q-1}$, the $\Zt$-action is the one which changes the
signs on the final $q$ coordinates.
Hence 
\[ X^\Zt=\{[r_1,\ldots,r_p,0\ldots,0] \mid r_i\in \R\} \amalg
         \{[0,\ldots,0,r_{p+1},\ldots,r_{p+q}] \mid r_i\in \R\};
\]
the first set is isomorphic to $\RP^{p-1}$, the second to $\RP^{q-1}$.

When $q\geq 1$, consider the map
\begin{eqnarray*}
  [(\R \oplus \R_{-})-0]/\pm 1 &\ra& X \\
   {[}r,s] &\assign& [r,0,\ldots,0,s,0,\ldots,0],
\end{eqnarray*}
where the $s$ is placed in the $q$th spot.  
It's easy to check that $[(\R\oplus \R_{-})-0]/\pm 1 \he S^{1,1}$, and
that on fixed sets the map induces a bijection on the sets of
path-components.  This is what we wanted.
\end{proof}

\vf

\appendix

\section{Symmetric products and their group completions}

The goal of this section is to show that $\SPi(S^{2n,n})$ is
equivariantly weakly equivalent to $\AG(S^{2n,n})$.  The proof is not
at all difficult, but requires a few lemmas.  I would like to thank
Gustavo Granja for an extremely helpful conversation about these results.
In this section $G$ is a fixed finite group.

\begin{defn}\mbox{}\par
\begin{enumerate}[(a)]
\item Let $\cC$ be a category with products and colimits, and let $M$
be an abelian monoid object in $\cC$.  The \dfn{group
completion} $M^+$ is the coequalizer of the maps
\[  \xymatrixrowsep{0.2pc}
    \xymatrix{ M\times M\times M \ar[r]<0.5ex>\ar[r]<-0.5ex> 
                      & M\times M \\
                      &{(a,b)}\\
               (a,b,c)\ar[ur]+/l 30pt/\ar[dr]\\
                      &{(a+c,b+c)}.
}
\]
\item If $K$ is a pointed simplicial set (or topological space),
define $\AG(K)=\Spi(K)^+$.
\end{enumerate}
\end{defn}

\begin{remark}
In the above generality it is not true that $M^+$ will be an abelian
group object in $\cC$, or even a monoid object (so the term `group
completion' is somewhat of a misnomer).  But this {\it is} the case
when $\cC=\Set$, and therefore also when $\cC=\sSet$.  It also
holds when $\cC=\Top$ and $M$ is `sufficiently nice'.
\end{remark}

It's easy to see that $\AG$ is a functor, so that if $K$ is a
simplicial $G$-set (or a $G$-space) then $\AG(K)$ also has a
$G$-action.  If $\tilde{\Z}[S]$ denote the free abelian group on the
pointed set $S$, where the basepoint is identified with the zero
element, one may check that $AG(K)$ is isomorphic to the simplicial
set $\tilde{\Z}[K]$.

\begin{prop}
\label{pr:simpspag}
Let $K$ be a pointed simplicial $G$-set with the property that $K^H$ is
path-connected for every subgroup $H\subseteq G$.  Then the natural
map $\SPi(K) \ra \AG(K)$ is an equivariant weak equivalence.
\end{prop}

\begin{proof}
If $M$ is a connected simplicial abelian monoid then 
\cite[Results Q1,Q2,Q4]{Q} show $M\ra M^+$ induces an equivalence on
integral homology.  Since both $M$ and $M^+$ are nilpotent spaces, the
map is actually a weak equivalence.  One can check, with only a little
trouble, that $AG(K)^H$ is isomorphic to $[\SPi(K)^H]^+$, and that
$\SPi(K)^H$ is path connected.  So Quillen's result implies that
$\SPi(K)^H \ra \AG(K)^H$ is a weak equivalence for every subgroup $H$.
This completes the proof.
\end{proof}

The next step is to transport this result from $G$-simplicial sets to
$G$-spaces.  We start with two simple lemmas:

\begin{lemma} Let $K$ be a pointed simplicial set.\par\noindent
\begin{enumerate}[(i)]
\item
There are natural homeomorphisms $\tre{SP^n K} \ra SP^n{\tre{K}}$, 
for $1\leq n\leq \infty$.
\item If $M$ is a simplicial abelian monoid, then there is a natural
homeomorphism $\tre{M^+} \ra \tre{M}^+$.
\item There is a natural homeomorphism $\tre{\AG(K)} \ra \AG(\tre{K})$. 
\end{enumerate}
\end{lemma}

\begin{proof}
The essential point is that $\tre{\blank}$ commutes with colimits
(being a left adjoint) and also finite products (making use of the
fact that our topological spaces are compactly-generated and weak
Hausdorff).  Since $\SP^n$ ($n<\infty$) is constructed by forming a
finite product and then taking a colimit---namely, passing to
$\Sigma_n$-orbits---realization will commute with $\SP^n$.  But then
realization will also commute with $\SPi$, as $\SPi$ is defined as a
colimit of the $\SP^n$'s.  This proves (i).

The proof of (ii) is in the same spirit: $M^+$ is defined as a
coequalizer of two products, so $\tre{\blank}$ will commute with this
construction.  Part (iii) is an immediate consequence of (i) and (ii).
\end{proof}

\begin{remark}
Again, since the above maps are all natural it follows that if $K$ is
a $G$-simplicial set then the maps are actually equivariant.    
\end{remark}

\begin{lemma} If $K$ is a $G$-simplicial set, the natural map
$\tre{K^H} \ra \tre{K}$ factors through the $H$-fixed set and gives a
homeomorphism $\tre{K^H} \ra \tre{K}^H$.
\end{lemma}

\begin{proof}
The fact that the map factors though $\tre{K}^H$ and that it is
injective are immediate.  So the content is that a cell of $\tre{K}$
which is fixed by $H$ must come from a simplex of $K^H$.  
But this is obvious.
\end{proof}

\begin{prop}
\label{pr:symm=ag}
Let $X$ be a $G$-space of the form $\tre{K}$ for some simplicial
$G$-set $K$.  If all the fixed sets $X^H$ are path-connected, then the
map $\SPi(X) \ra \AG(X)$ is an equivariant weak equivalence.
\end{prop}

\begin{proof}
What must be shown is that for any subgroup $H\subseteq G$ the map
$\SPi(X)^H \ra \AG(X)^H$ is a weak equivalence of topological spaces.
Now $X=\tre{K}$, and by the above lemmas we can commute the
realization past the $\SPi$, the $\AG$, and the fixed points.  So we
are left with showing that $\tre{\Spi(K)^H}\ra \tre{\AG(K)^H}$ is a
weak equivalence.  This was Proposition~\ref{pr:simpspag}.
\end{proof}

\begin{cor}
\label{co:symm=ag}
For any $V$ which contains a copy of the trivial representation,
the map of $\Zt$-spaces $\SPi(S^{V}) \ra \AG(S^{V})$ is an
equivariant weak equivalence.
\end{cor}

\begin{proof}
It suffices to show that $S^{V}$ is $\Zt$-homeomorphic to a space
of the form $\tre{K}$.
It's not hard to verify this for $V=\R$ and $V=\R_{-}$ 
by writing down an
explicit $K_1$ and $K_2$.  For a general $V=\R^p\oplus (\R_{-})^q$ we have 
\begin{eqnarray*}
 S^V &\iso& (S^{1,0} \Smash \ldots \Smash S^{1,0}) \Smash
(S^{1,1}\Smash \ldots \Smash S^{1,1}) \\
&\cong&
\tre{K_1}\Smash\ldots \Smash \tre{K_1} \Smash
\tre{K_2}\Smash\ldots \Smash \tre{K_2} 
\cong
\tre{
{K_1}\Smash\ldots \Smash {K_1} \Smash
{K_2}\Smash\ldots \Smash {K_2}}.
\end{eqnarray*}
\end{proof}

\vf 

\section{Computations of coefficient groups}

Here we give the proofs of Theorem~\ref{th:coeffs} and
Proposition~\ref{pr:HZandHZet}, which compute the coefficient rings of
$\HZ$ and $\Het$.  As remarked earlier, these computations are routine
among equivariant topologists---our only goal is to provide a
reference for the nonexpert.  

\medskip

\subsection{$H\Z$ computations}
For any pointed $\Zt$-spaces $X$ and $Y$ there is an
isomorphism $[\Zt_+ \Smash X,Y]_* \ra [X,Y]_*^e$ obtained by restricting
via the inclusion $\{0\} \inc \Zt$.  So for any equivariant spectrum
$E$ there are isomorphisms $E^{p,q}(\Zt) \ra E_{e}^{p}(\pt)$ where
$E_e$ is the  nonequivariant spectrum obtained by forgetting the group
actions.  If $E$ has a product, these isomorphisms give ring maps.  
It follows immediately that $H^{*,*}(\Zt)=\Z[u,u^{-1}]$ where $u$ has
degree $(0,1)$ (in effect, $u$ is just a placeholder for the second
index).  

For $H^{*,*}(\pt)$ we first recall that $H^{p,0}(\pt)$ is known by the
definition of Eilenberg-MacLane cohomology---it is $0$ when $p\neq 0$,
and $\Z$ when $p=0$.  For any space $X$ the groups $H^{p,0}(X)$ and
$H_{p,0}(X)$ are Bredon cohomology and homology with coefficients in
$\Z$, and in general one has $H^{p,0}(X)=H_{sing}^p(X/\Z_2;\Z)$ (but
the analog for homology is not quite true).
  
When $q>0$ we can now write 
\[ H^{p-q,-q}(\pt) \iso \tH^{p,0}(S^{q,q}) \iso
\tH^p_{sing}(S^{q,q}/\Z_2;\Z) \iso \tH^p_{sing}(\Sigma \RP^{q-1}).\]
(For the last isomorphism recall that $S^{q,q}$ is the suspension of
the sphere inside $\R^{q,q}$, which is the $(q-1)$-sphere with
antipodal action).  So when $q>0$ the groups $H^{*,-q}(\pt)$ are the
reduced cohomology of $\RP^{q-1}$, with a suitable shifting.  See the
picture in Theorem~\ref{th:coeffs}.

When $q>0$ we can also write
\[ H^{p+q,q}(\pt) \iso \tH^{p+q,q}(S^{0,0}) \iso
\tH_{-p-q,-q}(S^{0,0};\Z) \iso \tH_{-p,0}(S^{q,q}).\]
The second isomorphism uses equivariant Spanier-Whitehead duality.
Now, there is a cofiber sequence $S(\R^{q,q}) \inc D(\R^{q,q}) \ra
S^{q,q}$, where the first two terms are the sphere and disk in
$\R^{q,q}$.  The induced long exact sequence for $H_{*,0}$ shows that
$\tH_{a,0}(S^{q,q})\iso \tH_{a-1,0}(S(\R^{q,q}))$ when $a \neq 0,1$, and
that there is an exact sequence
\begin{equation}
\label{eq:long}
 0 \ra \tH_{1,0}(S^{q,q}) \ra H_{0,0}(S(\R^{q,q})) \ra \Z \ra
\tH_{0,0}(S^{q,q}) \ra 0.
\end{equation}
The space $S(\R^{q,q})$ has free action, and so
$\tH_{a-1,0}(S(\R^{q,q})) \iso
\tH_{a-1}^{sing}(S(\R^{q,q})/\Z_2)=\tH_{a-1}^{sing}(\RP^{q-1})$.
Hence $H^{p+q,q}(\pt)\iso \tH_{-p-1}^{sing}(\RP^{q-1})$ when $p\neq 0,1$.

The center map in (\ref{eq:long}) may be seen to coincide with the map
$H_{0,0}(\Zt) \ra H_{0,0}(\pt)$ induced by $\Zt \ra \pt$.  This is the
same as the transfer map $i_*$ in the Mackey functor $\Z$, which is
the $\times 2$ map $\Z \ra \Z$.  So $H^{q-1,q}(\pt)\iso
\tH_{1,0}(S^{q,q})=0$ and $H^{q,q}(\pt) \iso \tH_{0,0}(S^{q,q})\iso
\Zt$.  We have now seen that when $q>0$ the groups $H^{*,q}(\pt)$ are
the reduced singular homology groups of $\RP^{q-1}$, read downwards
from the group in degree $(q-1,q)$.  Again, the reader is referred to
the picture that goes with Theorem~\ref{th:coeffs}.

At this point we have computed the additive groups $H^{*,*}(\pt)$ and
$H^{*,*}(\Zt)$, so we turn our attention to the maps between them.
Consider the cofiber sequence $\Zt_+ \ra S^{0,0} \ra S^{1,1}$ and the
induced long exact sequence
\[  \cdots \la H^{0,2n-1}(\pt) \la H^{0,2n}(\Zt) \llla{i^*} H^{0,2n}(\pt) \la 
H^{-1,2n-1}(\pt) \la \cdots
\]
When $n\geq 0$ then $H^{0,2n-1}(\pt)=0=H^{-1,2n-1}(\pt)$, and so $i^*$ is an
isomorphism.  When $n<0$ we know enough to conclude that $\coker
i^*\iso \Zt$, so $i^*$ is
multiplication by $2$.  A $\Zt$-Mackey functor with both groups equal
to $\Z$ is completely determined by its restriction map $i^*$, so
we can deduce that $\underline{H}^{0,2n}=\Zm$ when $n\geq 0$ and
$\underline{H}^{0,-2n}=\Zm^{op}$ when $n<0$.

The map $i^*\colon H^{*,*}(\pt) \ra H^{*,*}(\Zt)$ is a map of rings,
and we know the target is $\Z[u,u^{-1}]$.  This allows us to determine
the subring $H^{0,2*}(\pt)$.  Also, the 
commutativity of the usual diagram
\[ \xymatrix{
AG(S^V) \Smash AG(S^W) \ar[r]\ar[d]^t & AG(S^V\Smash S^W)
 \ar[d]^{AG(t)} \\
 AG(S^W)\Smash AG(S^V) \ar[r] & AG(S^W\Smash S^V),}
\]  
shows that $H^{*,*}(X)$ is graded-commutative in a certain sense, for
any $X$.  For $X=\pt$ we know that the groups $H^{*,*}(\pt)$ are
either $\Zt$'s or else located in even degrees, so the ring is
commutative on-the-nose.  

It is not hard to see that $S^{1,1} \ra AG(S^{1,1})\he K(\Z(1),1)$ is
a weak equivalence (we know the homotopy groups of the target and its
fixed set, so this can be checked directly). 
Let $y$ denote the composite $S^{0,0} \inc S^{1,1} \ra K(\Z(1),1)$.
The cofiber sequence $S^{0,0} \inc S^{1,1} \ra \Zt_+ \Smash S^{1,0}$
gives us a long exact sequence on $H^{*,*}$ in which one of the maps
is multiplication by $y$.  Analysis of this long exact sequence lets
us determine all the multiplication-by-$y$'s shown in the diagram in
Theorem~\ref{th:coeffs}.  

At this point we have determined almost all of the ring structure on
$H^{*,*}(\pt)$.  If $\theta_n$ denotes the class in
$H^{0,-2n-1}(\pt)\iso\Zt$ and $x$ the generator of $H^{0,2}(\pt)\iso \Z$,
we have only to show that $x\cdot \theta_{n+1}=\theta_n$.  Let $E$ be
the spectrum defined by the cofiber sequence $\Sigma^{0,-2}H\Z \ra H\Z
\ra E$, where the first map denotes multiplication by $x$.  Using what
we have already proven, one computes that $E^{n,0}(\pt)=0$ if $n\neq
0$, $E^{0,0}(\pt)=\Zt$, and $E^{n,0}(\Zt)=0$ for all $n$.  So $E$ is
the Eilenberg-MacLane cohomology theory for the Mackey functor
$\underline{E}^{0,0}$, and the nature of this Mackey functor lets us
conclude that $E^{n,0}(X)\iso H_{sing}^n(X^{\Z_2};\Zt)$.  So when
$n>0$ we have $E^{0,-n}(\pt)\iso \tilde{E}^{n,0}(S^{n,n}) \iso
\tH^n_{sing}(S^0)=0$.  It follows that multiplication by $x$ gives an
isomorphism $H^{0,-n-2}(\pt) \ra H^{0,-n}(\pt)$ when $n\geq 2$.  This
completes the analysis of the ring structure on $H^{*,*}(\pt)$---all
products can be deduced from the ones we've computed together with
commutativity and degree considerations.

\subsection{$H\Z_{et}$ computations}
Recall that the spaces in the $\Omega$-spectrum for $\Het$ are
$K(\Z,V)^{E\Zt}$.  From this it follows that $\HZ \ra \Het$ is a
nonequivariant equivalence, and so $H^{*,*}(\Zt) \ra
H^{*,*}_{et}(\Zt)$ is an isomorphism of rings.

Remark~\ref{re:EMhfixed} gives the homotopy type of
$K(\Z(n),2n)^{h\Zt}$, and from this we immediately compute the groups
$H_{et}^{p,q}(\pt)$ where $p,q\geq 0$ and $p\leq 2q$.  The point is
that 
\[ H_{et}^{p,q}(\pt)=H_{et}^{2q-(2q-p),q}(\pt)=\tH_{et}^{2q,q}(S^{2q-p,0})=
[S^{2q-p,0},K(\Z(q),2q)^{h\Zt}].
\]

Using the cofiber sequence $S^{0,0} \inc S^{1,1} \ra \Zt_+\Smash
S^{1,0}$ now lets us deduce $H_{et}^{p,q}(\pt)$ in the two ranges
$(p\leq q)$ and $(p \geq 1)$.  In a moment we will show that for all
$n>0$ one has $H_{et}^{0,-2n}(\pt)=\Z$, $H_{et}^{0,-2n+1}(\pt)=0$, and
the restriction map $H_{et}^{0,-2n}(\pt) \ra H_{et}^{0,-2n}(\Zt)$ is
an isomorphism.  Using these facts, this same cofiber sequence
will show that $H_{et}^{p,q}(\pt)$ vanishes when both $p<0$ and $q<0$.

The above cofiber sequence induces an exact sequence 
\[ 0 = H_{et}^{1,-n+1}(\Zt) \la H_{et}^{1,-n+1} \llla{\cdot y} H_{et}^{0,-n} 
\la H_{et}^{0,-n+1}(\Zt) \llla{i^*} H_{et}^{0,-n+1}. 
\]
When $n=1$ we already know $i^*$ is the identity, and
$H_{et}^{1,0}=0$; so $H_{et}^{0,-1}=0$.  When $n=2$ we find the exact
sequence
\[ 0 \la \Zt \la H_{et}^{0,-2} \la \Z \la 0, 
\]
so $H_{et}^{0,-2}$ is either $\Z$ or $\Z \oplus \Zt$.
But we also know that
\begin{align*}
H_{et}^{0,-2} \iso \tH^{2,0}_{et}(S^{2,2})\iso
\tH^{2,0}(E\Zt_+\Smash S^{2,2})&\iso \tH^2_{sing}(E\Zt_+\Smash_{\Z_2}
S^{2,2})
\\
&\iso 
\tH^{2,0}_{sing}((S^1_a)_+\Smash_{\Z_2} S^{2,2}),
\end{align*}
where $S^1_a$ denotes the circle with antipodal action, and $S^1_a
\inc E\Z/2$ is the obvious inclusion.  The cofiber sequence
$(\Zt)_+ \inc (S^1_a)_+ \ra \Zt_+\Smash S^{1,0}$  gives a diagram of spaces
\[ \xymatrix{
\Zt_+ \Smash_{\Z_2} S^{2,2} \ar[r]\ar[d] &
   (S^1_a)_+ \Smash_{\Z_2} S^{2,2} \ar[r] \ar[dl] &
    (\Zt_+\Smash S^{1,0})\Smash_{\Zt} S^{2,2} \ar@{=}[r] & S^3 \\
S^{2,2}/\Z_2
}
\]
where the top row is a cofiber sequence and the two maps to the bottom
row squash $S^1_a$ (and $\Zt$) to a point. Applying
$H^2_{sing}$ to this diagram now gives
\[ \xymatrix{
H^3(S^3) & \Z \ar[l] & ? \ar[l] & 0\ar[l] \\
& \Z \ar[u]^{2}\ar[ur]
}
\]
where $?$ denotes $H^2_{sing}((S^1_a)_+\Smash_{\Z_2} S^{2,2})$.
We have so far determined that this group is either $\Z$ or
$\Z\oplus \Z/2$, and so the only possibility is $\Z$.  So we have
learned that $H^{0,-2}_{et}=\Z$, and the map $? \ra \Z$ in the diagram
must be an isomorphism.  The diagram now shows that $H^{0,-2} \ra
H^{0,-2}_{et}$ is multiplication by $2$, because this is the map $\Z
\ra ?$.  In the square
\[ \xymatrix{ H^{0,-2} \ar[r] \ar[d] & H^{0,-2}_{et} \ar[d] \\
  H^{0,-2}(\Zt) \ar[r] & H^{0,-2}_{et}(\Zt) 
}
\]
we know all the groups are $\Z$, the top and left maps are
multiplication by $2$, and the bottom map is an isomorphism; so the
right vertical map is also an isomorphism.

Using an induction, the above arguments actually show that
$H_{et}^{0,-2n+1}=0$, $H_{et}^{0,-2n}=\Z$, and $H^{0,-2n}_{et} \ra
H^{0,-2n}_{et}$ is an isomorphism, for all $n\geq 1$.  
This completes our determination of the groups $H^{*,*}_{et}$, and of
the Mackey functors $\underline{H}_{et}^{0,2n}$.  

The ring structure on $H^{0,*}_{et}(\pt)$ can now be determined by
comparing with the known structure of the rings $H^{0,*}(\pt)$ and
$H^{0,*}(\Zt)$ (the latter is also $H^{0,*}_{et}(\Zt)$).  The
multiplication-by-$y$'s are deduced from the long exact sequences
induced by $S^{0,0}\inc S^{1,1} \ra \Zt_+\Smash S^{1,1}$, just as for
$H^{*,*}$.  This completes the proof of
Proposition~\ref{pr:HZandHZet}.

\bibliographystyle{amsalpha}

\begin{thebibliography}{JTTW}

\bibitem[AM]{AM} S. Araki and M. Murayama, \emph{$\tau$-cohomology
theories}, Japan. J. Math {\bf 4} (1978), no. 2, pp. 361--416.

\bibitem[Ar]{Ar} S. Araki, \emph{Forgetful spectral sequences}, Osaka
J. Math. {\bf 16} (1979), pp. 173--199.

\bibitem[AH]{AH} M. Atiyah and F. Hirzebruch, \emph{Vector bundles and
homogeneous spaces}, Proc. Symp. Pure Math. {\bf 3}, pp. 7--38. 

\bibitem[At]{At} M. Atiyah, \emph{$K$-theory and reality}, Quart. J.
Math. Oxford (2) {\bf 17} (1966), pp. 367--386.

\bibitem[Be]{Be} A. Beilinson, \emph{Height pairing between algebraic
cycles}, in \emph{K-theory, Arithmetic and Geometry}, Lect. Notes in
Math. {\bf 1289}, Springer (1987), pp. 1--26.




\bibitem[BL]{BL} S. Bloch and S. Lichtenbaum, \emph{A spectral
sequence for motivic cohomology}, preprint, 1995.  Available at
http://www.math.uiuc.edu/K-theory/0062.

\bibitem[Bd]{Bd} M. Boardman, \emph{Conditionally convergent spectral
sequences}, in \emph{Homotopy invariant algebraic structures
(Baltimore, MD, 1998)}, 49--84, Contemp. Math. {\bf 239},
Amer. Math. Soc., Providence, RI, 1999.


\bibitem[BK]{BK} A. Bousfield and D. Kan, \emph{Homotopy limits,
completions, and localizations}, Lect. Notes in Math. {\bf 304},
Springer-Verlag, 1972.

\bibitem[BG]{BG} K. Brown and S. Gersten, \emph{Algebraic $K$-theory
and generalized sheaf cohomology}, Lecture Notes in Math., vol. 341,
Springer, 1973, pp. 266--292.

\bibitem[dS]{dS} P. dos Santos, \emph{A note on the equivariant
Dold-Thom theorem}, to appear in J. Pure. Appl. Alg.

\bibitem[D1]{D1} D. Dugger, \emph{A Postnikov tower for algebraic
$K$-theory}, MIT PhD thesis, 1999.

\bibitem[D2]{D2} D. Dugger, \emph{Multiplicative structures on
homotopy spectral sequences I, II}, preprint, 2003.

\bibitem[DI]{DI} D. Dugger and D. Isaksen, \emph{Hypercovers in
topology}, to appear in Math. Zeit.

\bibitem[F]{F} E. Dror Farjoun, \emph{Cellular Spaces, Null Spaces,
and Homotopy Localization}.  Lect. Notes in Math. {\bf 1622},
Springer-Verlag, 1996.

\bibitem[FS]{FS} E. Friedlander and A. Suslin, \emph{The spectral
sequence relating motivic cohomology to algebraic $K$-theory},
Ann. Sci. \'Ecole Norm. Sup. (4) {\bf 35} (2002), no. 6, 773--875.

\bibitem[Gr]{Gr} D. Grayson, \emph{Weight filtrations in algebraic
$K$-theory}, in \emph{Motives}, Proc. Symp. Pure Math. {\bf 55} (1994), pp.
207--244.  


\bibitem[H]{H} P. Hirschhorn, \emph{Model Categories and Their
Localizations}, Mathematical Surveys and Monographs, vol. 99,
Amer. Math. Soc., 2003.

\bibitem[HK1]{HK1} P. Hu, I. Kriz, \emph{Real-oriented homotopy theory
and an analogue of the Adams-Novikov spectral sequence}, Topology {\bf
40} (2001), no. 2, 317--399.

\bibitem[HK2]{HK2} P. Hu, I. Kriz, \emph{Some remarks on Real and
algebraic cobordism}, $K$-theory {\bf 22} (2001), no. 4, 335--366.


\bibitem[JTTW]{JTTW} I. James, E. Thomas, H. Toda, G. W. Whitehead,
\emph{On the symmetric square of a sphere}, J. Math. Mech. {\bf 12}
(1963), pp. 771--776.

\bibitem[J]{J} J.F. Jardine, \emph{Motivic symmetric spectra},
Doc. Math. {\bf 5} (2000), 445--553.

\bibitem[K]{K} M. Karoubi, \emph{A descent theorem in topological
$K$-theory},  K-theory {\bf 24} (2001), no. 2, 109--114.

\bibitem[LLM]{LLM} H. B. Lawson, P. Lima-Filho, M. Michelsohn,
\emph{Algebraic cycles and the classical groups Part I, Real cycles},
Topology 42 (2003), no. 2, 467--506.
 
\bibitem[Lw1]{Lw1} L. G. Lewis, \emph{The $RO(G)$-graded equivariant
ordinary cohomology of complex projective spaces with linear $\Z/p$
actions}, appearing in \emph{Algebraic topology and transformation
groups}, ed. tom Dieck, Lecture Notes in Math., vol. 1361, Springer,
1988.
 
\bibitem[Lw2]{Lw2} L. G. Lewis, \emph{The equivariant Hurewicz map},
Tran. Amer. Math. Soc., {\bf 329} (1992), no. 2, pp. 433--472.

\bibitem[Lw3]{Lw3} L. G. Lewis, \emph{Equivariant Eilenberg-MacLane spaces
and the equivariant Seifert-van Kampen suspension theorems},
Topology Appl., {\bf 48} (1992), no. 1, pp. 25--61.

\bibitem[LMS]{LMS} L. G. Lewis, J. P. May, and M. Steinberger, {\emph
Equivariant stable homotopy theory}, Lecture Notes in Math.,
vol. 1213, Springer, 1986.

\bibitem[LF]{LF} P. Lima-Filho, \emph{On the equivariant homotopy of
free abelian groups on $G$-spaces and $G$-spectra}, Math. Zeit.
{\bf 224} (1997), pp. 567--601.


\bibitem[M]{M} J.P. May, et al, \emph{Equivariant Homotopy and
Cohomology Theory}, CBMS volume 91, American Mathematical Society,
1996.

\bibitem[MV]{MV} F. Morel and V. Voevodsky, \emph{$\A^1$-homotopy
theory for schemes}, Inst. Hautes \'Etudes Sci. Publ. Math., No. 90
(2001), 45--143.

 
\bibitem[Q]{Q} D. Quillen, \emph{On the group completion of a
simplicial monoid}, appearing as an appendix to \emph{Filtrations on
the homology of algebraic cycles} (Friedlander and Mazur), Mem. Amer.
Math. Soc. (1994), no. 529.

\bibitem[V1]{V1} V. Voevodsky, \emph{The Milnor conjecture}, preprint,
1996.  Available from the $K$-theory archive at
http://www.math.uiuc.edu/K-theory/0170.

\bibitem[V2]{V2} V. Voevodsky, \emph{Open problems in motivic stable
homotopy theory, I\/}, preprint, 2000.  Available at 
http://www.math.uiuc.edu/K-theory/0392. 

\bibitem[V3]{V3} V. Voevodksy, \emph{A possible new approach to the
motivic spectral sequence}, appearing in \emph{Recent progress in
homotopy theory (Baltimore, MD, 2000)}, 371--379, Contemp. Math. 293,
Amer. Math. Soc., Providence, RI, 2002.


\bibitem[V4]{V4} V. Voevodksy, \emph{On the zero-slice of the sphere
spectrum}, preprint, 2002.  Available at
http://www.math.uiuc.edu/K-theory/0612.

\bibitem[W]{W} S. Waner, \emph{$G-CW(V)$ complexes}.  Unpublished manuscript.

\end{thebibliography}

\end{document}